\documentclass[10pt]{amsart}
\usepackage{amsmath,amscd}

\title{Indecomposable Higher Chow cycles on Jacobians}

\author{Alberto Collino}
\address{Dipartimento di Matematica \\ Universit\`a di Torino \\ Via Carlo Alberto 10 \\ 10123 Torino \\ Italy}

\email{collino@dm.unito.it}
\author{Najmuddin Fakhruddin}
\email{nfakh@bom2.vsnl.net.in}

\thanks{First author partially supported partially supported by
Progetto nazionale MURST 1997 "Algebra commutativa, geometria 
algebrica e applicazioni", UE Research Training Network EAGER}


\newtheorem{thm}{Theorem}[section]
\newtheorem{lem}{Lemma}[section]
\newtheorem{prop}{Proposition}[section]
\newtheorem{cor}{Corollary}[section]

\newtheorem{hyp}{Hypothesis}[section]

\theoremstyle{remark}
\newtheorem*{rem}{Remark}

\newtheorem*{conv}{Conventions}
\newtheorem*{ack}{Acknowledgements}

\theoremstyle{definition}
\newtheorem*{defn}{Definition}
\newtheorem*{assumption}{Assumption}

\numberwithin{equation}{section}

\begin{document}
\addtocounter{section}{-1}

\begin{abstract}
We construct some natural cycles with trivial
regulator in the higher Chow groups
of Jacobians. For hyperelliptic curves we use a criterion due to
J.~Lewis to prove that the cycles we construct are indecomposable,
and then use a specialization
argument to prove indecomposability for more general curves.
\end{abstract}

\maketitle

\section{Introduction}

The aim of this paper is to construct some natural cycles
in the higher Chow groups of Jacobians of smooth projective
curves. Most of the paper is devoted to the first higher
Chow groups $CH^k(J(C),1)$,
especially the case $k = g = genus(C)$, but in the
last section we also construct elements of  $CH^k(J(C),n)$, $n>1$,
for curves of low genus. 

For a smooth projective variety $X$,
the subgroup of  $CH^k(X,1)$ of decomposable cycles,
 $CH^k_{\mathrm{dec}}(X,1)$, is 
defined to be the image of $CH^{k-1}(X) \otimes \mathbb{C}^{*}$,
where we use  the isomorphism $CH^1(X,1) \cong\mathbb{C}^{*}$. 
We let $CH^k_{\mathrm{ind}}(X,1) :=
CH^k(X,1)\big/CH^k_{\mathrm{dec}}(X,1)$ be the quotient group
of indecomposable cycles. In \cite{collino-jag}, a natural
element $K$ was constructed in $CH^g(J(C),1)$,  $C$ 
hyperelliptic, and using the regulator map
$K$ was shown to give a non-torsion element
of $CH^g_{\mathrm{ind}}(J(C),1) $ when $C$ is generic hyperelliptic.
The Pontryagin product of $K$ with zero cycles of degree zero,
gives elements of  $CH^g(J(C),1)$ which lie in the kernel of
the regulator map. Our first result, Theorem \ref{thm:hypmain},
shows that such cycles give uncountably many elements of
$CH^g_{\mathrm{ind}}(J(C),1)_{\mathbb{Q}}$ if $g \geq 3$. 
We prove a more precise statement using the decomposition
of the higher Chow groups due to Beauville \cite{be2} and 
Deninger-Murre \cite{dm}. Applying the motivic hard
Lefschetz theorem of K\"unnemann \cite{kunnemann-lefs}
we obtain uncountably many elements of
 $CH^k_{\mathrm{ind}}(J(C),1)_{\mathbb{Q}}$, $ 3 \leq k \leq g$,
lying in the
kernel of the regulator map. The main
technical tool used in the proof is a Hodge theoretic criterion due to J.~Lewis
\cite{lewis-crelle}. This has been used earlier by Gordon and Lewis
\cite{gordon-lewis-jag} to construct indecomposable cycles
with similar properties in products of generic elliptic curves.

Our main goal is to construct indecomposable cycles in 
$CH^g(J(C),1)$ for more general curves.
A result of Nori \cite[7.5]{nori-invent} implies that upto torsion
the regulator image of
$CH^g(J(C),1)$ is the same as that of $CH^g_{\mathrm{dec}}(J(C),1)$
for a generic curve of genus $g=3$, and
it seems likely that this should also be true for higher $g$.
This makes it difficult to use Lewis' criterion to prove
indecomposability; instead we employ a specialization
argument.
To this end, we first prove that the specialization of
a decomposable cycle is decomposable (Theorem \ref{thm:decomp}).
The difficulty here is that a cycle on a family of smooth
projective varieties which restricts to a decomposable cycle
on a generic fibre need not be decomposable on the entire
family. We circumvent this problem by considering an
auxiliary family which is constructed by gluing a 
product family along a special fibre.

Given a divisor $D = (a_1 + a_2) - (b_1 + b_2)$ on a curve $C$
of genus $g \geq 2$ with $2[D] = 0$ in $J(C)$, we construct 
an element $Z_D$ of $CH^g(J(C),1)$,  called the $4$-configuration since it is
supported on $4$ copies of $C$
in $J(C)$. When $C$ and $D$ are generic and $g \geq 3$, we show
(Theorem \ref{thm:indec4} that this element is indecomposable: By 
using a suitable Hurwitz
scheme, we show that we may specialize $C$ to  a hyperelliptic
curve  $C'$ in such a way that the $4$-configuration specializes
to $K - K_t$, with $t$ a generic point of $C'$. Indecomposability
follows by applying Theorem \ref{thm:hypmain} and \ref{thm:decomp}.
We then deduce that $CH^g_{\mathrm{ind}}(J(C),1)_{\mathbb{Q}}$ is
uncountable for generic curves of genus $3$ and $4$.

To construct elements of $CH^g(J(C),1)$ for arbitrary curves,
we consider divisors $D = \sum_{i = 1}^n a_i - \sum_{i=1}^n b_i$
with $2[D] = 0$ in $ J(C)$ and $n>0$. Associated to such a divisor
there is a natural subspace of $CH^{n+1}(C^{n+1},1)_{\mathbb{Q}}$
and we describe a method, generalizing that used for the 
$4$-configuration, which we believe should enable one to prove
indecomposability of the general such element with $C$ and $D$ also
generic and $g \geq n+1$. Unfortunately, due to the combinatorial
complexity of the cycles involved, we have not been able to
complete the proof. However, we have checked indecomposability using
a simple computer program for  $n \leq 6$. In particular,
we see that 
$CH^g_{\mathrm{ind}}(J(C),1)_{\mathbb{Q}}$ is uncountable for
a generic curve with $3 \leq genus(C) \leq 12$.

We conclude the paper with a construction of elements
in $CH^3(J(C), 4-g)$, with $g \leq 2$. These may be viewed
as successive degenerations of the $4$-configuration on a genus $3$
curve, as the curves acquire nodes. We expect, but do not prove,
that these elements are indecomposable in a sense stronger 
than that of Lewis. We do prove however,
using Lewis's criterion, that for $B$ a bielliptic genus $2$
curve $CH^3(J(B),2) \big /Im(K_2(\mathbb{C}) \otimes CH^1(J(C)))$
is non-trivial modulo torsion.

\begin{conv}
All varieties will be over the complex numbers $\mathbb{C}$
and all points will be closed points. We shall
say that a condition holds for a \emph{general} point
of a variety if it holds for all points in
a Zariski open subset and it holds for a \emph{generic} point if
it holds for all points outside a countable union of proper subvarieties.

We denote by $*$ the Pontryagin product
on the higher Chow groups of an abelian variety.
\end{conv}

\section{ Lewis' conditions, Pontryagin products and hyperelliptic Jacobians}
\label{sec:hyp}

The first higher Chow group $ CH^{k}(X,1) \simeq H^{k-1}(X,\mathcal
K_{k})$ of a non singular variety $X$ is generated by higher cycles of
the form $ Z = \sum_i Z_i \otimes f_i$, where the $Z_i$ are
irreducible subvarieties of codimension $(k-1)$ and the rational
functions $f_i \in k(Z_i)^\times$ obey the rule $\sum_i div(f_i) =0$
as a cycle on $X$.  Consider the subgroup of \emph{decomposable}
cycles
\begin{equation}
CH_{\mathrm{dec}}^k(X,1) := Im \big\{ CH^1(X,1) {\otimes} CH^{k-1}(X)
\longrightarrow CH^k(X,1) \big\} \ , \label{def:decomp}
\end{equation}
and the related quotient of \emph{indecomposable} cycles
\begin{equation*}
CH_{\mathrm{ind}}^k(X,1) := CH^k(X,1)\big/
CH_{\mathrm{dec}}^k(X,1) \ .
\end{equation*}
\noindent Recall that if $X$ is projective,
then $CH^1(X,1) = \mathbb{C}^{*}$.

Let $I$ be the zero cycles of degree zero on an abelian variety $A$.
Bloch has shown that $I^{*n}$ is
non-zero for $1 \leq n \leq g$
whereas $I^{*(g+1)}$ is always zero, $g = dim(A)$. Given $Z$ an
indecomposable element in $CH^k(A,1)$, a natural question to
ask is whether $I^{*n}*Z $ contains indecomposable cycles for $n\geq 1$.  We
remark that $I^{*n}*Z $ is in the kernel of the regulator map
for any $Z$, since
translation acts trivially on the Deligne cohomology
$ H_{\mathcal{D}}^{2k-1}(X,\mathbb{Z}(k))$.

Using as a basic tool a condition of Hodge type due to Lewis,
for $Z$ a real regulator indecomposable element of $CH^g(J(C),1)$
we give
a criterion in terms of the primitive cohomology of $J(C)$ for
$I^{*n}*Z $ to contain indecomposable cycles.
A significant instance of this situation is the case of generic
hyperelliptic Jacobians, where we can use as $Z$ the basic cycle in
$CH_{\mathrm{ind}}^g(J(C),1)$ found in \cite{collino-jag}.
We show that  $I^{*n}*Z $ contains indecomposable cycles for
$1 \leq n \leq g-2$, whereas  all elements of $I^{*(g-1)}*Z $
are decomposable (with $\mathbb{Q}$ coefficients).

\subsection{Preliminaries} 

We recall some notation and definitions and then we state a theorem of
Lewis. Our aim is to have a concrete reference at hand, for
more details the reader should consult either the original paper
\cite{lewis-crelle} or the survey \cite{gordon-lewis-banff}.

For $X$ projective and nonsingular,
\begin{equation*}
  \frac{H^{i-1}(X,\mathbb C)}{F^jH^{i-1}(X,\mathbb C) + H^{i-1}(X,\mathbb
    R(j))} \simeq
  \frac{H^{i-1}(X,\mathbb R(j-1))}{\pi_{j-1}(F^jH^{i-1}(X,\mathbb C))} 
\end{equation*}
and therefore one has the identification
\begin{align*}
  H_{\mathcal D}^{2k-1}(X,\mathbb R(k)) &\simeq H^{2k-2}(X,\mathbb
  R(k-1)) \cap F^{k-1}H^{2k-2}(X,\mathbb C) \\ \vspace{1\jot} &=:
  H^{k-1,k-1}(X,\mathbb R(k-1)) \, .
\end{align*}

According to Beilinson \cite{beilinson-ahc}, the real regulator image
of a cycle $ Z= \sum Z_i \otimes f_i$ in $ CH^k(X,1)$ is the
element
\begin{equation*}
  R_{k,1}( Z) \in  H_{\mathcal D}^{2k-1}(X,\mathbb R(k)) \simeq
  H^{k-1,k-1}(X,\mathbb R(k-1))\, ,
\end{equation*}
determined by the class of the current
\begin{equation*}
 R_{k,1}( Z) : \omega \longmapsto (2\pi \sqrt{-1})^{k-1-d} \sum_i
\int_{Z_i -
Z_i^{\mathrm{sing}}} \omega \,\log|f_i|\,  .
\end{equation*}

\begin{defn}
  A higher Chow cycle $ Z \in CH^k(X,1) $ is said to be
  (real) {\it regulator indecomposable} if there exists a
  differential form
  \begin{equation*}
    \omega \in (Hdg^{k-1}(X) \otimes \mathbb R )^{\perp} \subset
    H^{d-k+1,d-k+1}(X,\mathbb R(d-k+1))
  \end{equation*}
  such that the pairing $[R_{k,1}(Z) , \omega] \ne 0$.
\end{defn}

\begin{lem}  If $ Z \in   CH^k(X,1)$ is regulator
  indecomposable then $ Z$ is indecomposable.
  \end{lem}
  
Recall that the {\it coniveau filtration} on $H^i(X,\mathbb Q)$ is
\begin{equation*}
  N^jH^i(X,\mathbb Q) := \ker \big(\,H^i(X,\mathbb Q) \longrightarrow
  lim_{codim{_X} Y\ge j} H^i(X-Y,\mathbb Q) \,\big)\, ,
\end{equation*}
where the direct limit is over closed subvarieties $Y \subset X$. The
complex subspace generated by the Hodge projected image of the
coniveau filtration is
\begin{equation*}H_N^{k-l,k-m}(X) :=Im (\, N^{k-l}H^{2k-l-m}
  (X,\mathbb Q)\otimes \mathbb C \longrightarrow
  H^{k-l,k-m}(X) \,\big) \, .
\end{equation*} 
\noindent Lewis constructs certain complex subspaces
\begin{equation*}
  H^{\{k,l,m\}}(X) \subseteq  H^{k-l,k-m}(X) \, ,
\end{equation*}
such that for $m=0$ one has
\begin{equation}
  H^{\{k,l,0\}}(X) \subseteq  H_N^{k-l,k}(X) \, .
  \label{eq:coniveau}
\end{equation}
\noindent The spaces $H^{\{k,l,m\}}(X)$ are obtained by a process of K\"unneth
projection of the Hodge components of the real regulator classes of
the elements in $CH^k(X \times S,m)$, with varying smooth projective
varieties $S$.

We give an abridged version of the main result from
\cite{lewis-crelle}, since we only need it for
$CH_{\mathrm{ind}}^k(X,1)$.
\begin{thm}[Lewis]
  Let $X$ be a non singular projective variety:
\begin{enumerate}
\item $H^{\{k-1,l-1,0\}}(X) \subset H^{\{k,l,1\}}(X)$.  \smallskip
\item If $H^{\{k,l,1\}}(X) / H^{\{k-1,l-1,0\}} (X) \ne 0$ for some $l$,
  $2 \le l \le k$, then $CH_{\mathrm{ind}}^k(X,1)_{\mathbb{Q}}$ is
  uncountable.
\end{enumerate}
\end{thm}
Note that Lewis' proof shows that for a cycle $\xi \in CH^k(X \times
S,1)$ which gives rise to a nonzero element of $H^{\{k,l,1\}}(X) /
H^{\{k-1,l-1,0\}} (X)$, the cycles $\xi_s$, $s \in S$, form an
uncountable subset of $CH^k_{\mathrm{ind}}(X,1)_{\mathbb{Q}}$.

Our results will show that (ii) holds for a generic hyperelliptic
Jacobian of genus $g \geq 3$, for $k=g$ and all $l$ such that 
$2 \leq l \leq g-1$.

\subsection{Lewis' condition holds for the Pontryagin families $Z(m)$} 
\label{subsec:zm}

Let $C$ be a smooth projective curve of genus $g$ and $J = J(C)$ be
its Jacobian.  Let $Z \in CH^ { g} (J(C), 1)$ be a regulator
indecomposable cycle.  Given $m$ points $p_i$ on $C$, we construct
higher cycles $Z(m)$ in $CH^ { g} (C^{m}\times J(C), 1)$, so that the
fibre over $ (t_1, \dots , t_m) $ is $ Z * ( [t_1 - p_1] - e) * \dots
*( [t_m - p_m] -e) $, where $e$ is the origin in $J(C)$. Set $Z(1) :=
b_*( C \times Z ) - C \times Z $, where $ b : C\times J(C) \to C\times J(C)
$ is the twisted isomorphism $b (t,x) = (t,[t - p_1] +x )$.  By
iteration this gives $Z(m)$ on $C^{m}\times J(C) = C \times (C^{m-1}
\times J(C)) $, where the twisted map is now $ b( (t_m,\dots,t_1,x) )=
(t_m,\dots,t_1,[t_m -p_m] +x )$.

\begin{prop} Let $C$ be a smooth projective curve of genus $g$ and $Z$ a
  regulator indecomposable cycle in $CH^{g}(J(C), 1)$. Assume that the
  primitive cohomology $P^i(J(C))$ is an irreducible sub-Hodge structure
  of $H^i(J(C),\mathbb Q)$ for some $m \leq g-2$ and all $i \leq m+2$.
  Then the real regulator of $Z(m)$ gives rise to a non-zero element of
  $H^{\{g,m+1,1\}}(J(C)) / H^{\{g-1,l,0\}}(J(C))$, and hence the 
  restriction of $Z(m)$ to the fibre over a generic point of $C^m$
  is indecomposable.
  \label{prop:zm}
\end{prop}

If $ \omega_1^{J} , \dots, \omega_g^{J}$ is a basis for $H^0(J(C),\Omega
^{1}_J(C))$ such that the restriction $\zeta_i := \omega_{i \, |C}^{J}$,
$1 \leq i \leq g$, is an orthonormal frame for $H^0(C,\Omega ^{1}_C)$,
then the class of the divisor $\Theta$ is determined by the form $
\theta_J(C) = (i/2) \sum_{j=1}^{g} \omega_j^{J}
\wedge\bar{\omega}_j^{J}$.  We define $ \tau_C = $ $\omega_1^{J}
\wedge\bar{\omega}_1^J -\omega_2^{J} \wedge\bar{\omega}_2^J$.

\begin{lem}  
  Under the assumptions of Proposition \ref{prop:zm}, if Z is real
  regulator indecomposable then there is a basis as above with $[R(Z),
  \tau] \neq 0$.
\end{lem}
\begin{proof}
  The inner product on $H^0(C,\Omega^{1}_C)$ allows us to view it as a
  representation $V$ of the unitary group $U(g)$. $H^{1,1}(J(C))$ is then
  also a representation of $U(g)$, and is isomorphic to a twist of
  $V \otimes V^*$ 
  by a 1-dimensional representation.  Hence by Pieri's formula we see
  that it decomposes as a direct sum of two irreducible
  representations $T$ and $U$, $T$ being the subspace corresponding to
  the class of the $\Theta$ divisor.
  
  Now consider the subspace $W$ of $H^{1,1}(J(C))$ spanned by all
  possible $\tau$'s as above. Since the unitary group preserves the
  inner product, it follows that $W$ is a subrepresentation of
  $H^{1,1}(J(C))$. Clearly $W$ is not equal to $T$, hence it must contain
  $U$. If there were no $\tau$'s with a non-zero pairing with $Z$ then
  the pairing would be zero on all of $W$, contradicting the
  assumption of real regulator indecomposable. (Note that our
  assumptions imply that the space spanned by the rational Hodge
  classes is $T$).

\end{proof}

Let $\alpha_m =\bar \zeta_{m+2}\dots \bar \zeta_3 \wedge \tau \wedge
\omega_3^{J} \wedge \dots \omega_{m+2}^{J} $. 

\begin{lem}  \label{lem:alpha}
  $[ R(Z(m)),\alpha_m ] = (-1) ^{m} [ R(Z),\tau] $
\end{lem}
\begin{proof}
  By iteration the proof is the same for all $m \geq 1$. Say $m=1$,
  then we have $b^{*} (\omega_s^{J}) = \omega_s^{J} + \zeta_s $ and $b^{*}
  (\zeta_s )= \zeta_s $, thus 
  \[
  b^* (\bar \zeta_3 \wedge \tau \wedge \omega_3^{J} ) = \bar
  \zeta_3 \wedge b^* ( \tau \wedge \omega_3^{J} ) = -(\zeta_3
  \wedge\bar{\zeta}_3) \wedge \tau + \sum_i \phi_i^{C} \wedge
  \psi_i^{J} \, .
  \]
  Here $\phi_i^{C}$ is a form from C and $\psi_i^{J}$ is
  from $J(C)$, and it is either $\phi_i^{C} = \bar \zeta_3 $ or else
  $\phi_i^{C} = \bar \zeta_3 \wedge \zeta_l , l\neq 3 $, and therefore
  it is of volume $0$, because of the orthogonality assumption. We
  have then : $[ R(Z(1) ), \bar \zeta_3 \wedge \tau \wedge
  \omega_3^{J}] = $ $[ R( C\times Z), b^{*} (\bar \zeta_3 \wedge \tau
  \wedge \omega_3^{J} ) ] - [ R( C\times Z), \bar \zeta_3 \wedge \tau
  \wedge \omega_3^{J} ] = $ $ [ R( C\times Z), b^{*} (\bar \zeta_3
  \wedge \tau \wedge \omega_3^{J} ) ] = [-R(Z),\tau] $, because $ [ R(
  C\times Z), \phi^{C} \wedge \psi^{J} ]= $ $[ R( Z), \psi^{J}]
  \int_{C} \phi^{C} $.
\end{proof}

\begin{proof}[Proof of Proposition \ref{prop:zm}]
  The previous lemma along with equation \eqref{eq:coniveau} implies
  that $R(Z(m))$ gives a non-zero element of $H^{\{g,m+1,1\}}(J(C)) /
  H^{\{g-1,m,0\}}(J(C))$ if $\tau \wedge \omega_3^{J} \wedge \dots
  \omega_{m+2}^{J}$ is orthogonal to $N^{g-m-1}H^{2g-m-2}(J(C))$.  The
  assumptions on the primitive cohomology imply that
  $N^{g-m-1}H^{2g-m-2}(J(C)) = \Theta^{g-m-1}H^m(J(C))$.  Since $\tau \wedge
  \omega_3^{J} \wedge \dots \omega_{m+2}^{J} \cdot \Theta^{g-m-1} =0$,
  it follows that the condition is indeed satisfied.
\end{proof}

The next proposition shows that the hypothesis on the primitive
cohomology of  Proposition \ref{prop:zm} holds
for the Jacobian of a generic hyperelliptic curve.
 The reader may refer to
\cite{mumford-hodge} for the definition and basic properties of the
Hodge group.

\begin{prop}
  Let $C$ be a generic hyperelliptic curve of genus $g$. Then the
  Hodge group of $J(C)$ is isomorphic to $Sp(2g, \mathbb Q)$,
  hence $P^i(J(C))$
  is an irreducible Hodge structure for  $0 \leq i \leq g$.
  \label{prop:hodgegroup}
\end{prop}

\begin{proof}
  We shall use induction on $g$, the result for $g =1$ and $2$ being
  well known.  Assume that $g \geq 3$ and the result is known for
  smaller $g$.  We degenerate $C$ to a stable curve $C_o$ with $3$
  smooth irreducible components $C_1, C_2$ and $C_3$, with $C_1$ and
  $C_3$ of genus $1$, and $C_2$ of genus $g-2$ intersecting each of
  $C_1$ and $C_3$ transversally in a single Weierstrass point.  By
  choosing a path in a suitable parameter space, we can identify $V :=
  H^1(C, \mathbb Q)$ as a symplectic vector space with $ V_1 \oplus
  V_2 \oplus V_3$, where $V_i = H^1(C_i, \mathbb Q)$, $i=1,2,3$. Let
  $D_1$ and $D_3$ be generic hyperelliptic curves of genus $g-1$, with
  $D_i$ specializing to $C_i \cup C_2$, $i=1,2$. Let $C' = D_1 \cup
  C_3$ and $C'' = C_1 \cup D_3$, the two components of each curve
  intersecting transversally in a single point.  Again, by choosing
  paths we may identify $H^1(C', \mathbb Q)$ and $H^1(C'', \mathbb Q)$
  with $ V_1 \oplus V_2 \oplus V_3$ in such a way that $H^1(D_i,
  \mathbb Q)$ is identified with $V_i \oplus V_2$.
  
  Now consider the family of Jacobians. Since $C$ is generic, $G$, the
  Hodge group of $J(C)$, contains the Hodge groups $G'$ and $G''$ of
  $J(C')$ and $J(C'')$ respectively. Using induction and the above
  identifications, we see that $G$ contains both $Sp(V_1 \oplus V_2)$
  and $Sp(V_2 \oplus V_3)$.  One easily checks, by an explicit
  computation using Lie algebras, that the smallest subgroup of
  $GL(V)$ containing both these two subgroups is $Sp(V)$. The
  Hodge group is always contained in
  $Sp(V)$ so $G$ must equal $Sp(V)$.

  From the representation theory of symplectic groups it follows that
   $P^i(J(C))$ is an irreducible representation of the
  Hodge group, $0 \leq i \leq g$. Since the sub-Hodge structures of
  $H^i$ are precisely the subrepresentations of the Hodge group, it
  follows that the $P^i(J(C))$'s are also irreducible as Hodge structures.
\end{proof}
 
\subsection{On the  real regulator image of the basic hyperelliptic cycle} 
\label{subsec:hypreg}
Let $f : C \to \mathbb P^{1}$ be the double cover associated with a
hyperelliptic curve. We fix two ramification points $w_1$ and $w_2$ on
$C$ and choose a standard parameterization on $\mathbb P^{1}$ so that
$ f(w_1)=0$ and $ f(w_2) = \infty $.  The points $w_1$ and $w_2$ are
the distinguished Weierstrass points, and $ \epsilon := [{w_1}-{w_2}]$
is the associated element of order two in $Pic^0(C) $.
 
It is convenient to identify $J(C)=Pic^0(C) $ with $Pic^1(C)$ by
adding $w_1$. We embed $C$ in the natural way in $Pic^1(C)$, and for $t
\in Pic^1(C)$ we let $C_t$ be the translate of $ C $ by $[t-w_1]$.
Now consider $ W_1 := C = C_{w_1}$ and $ W_2 :=C_{w_2}$, the $
\epsilon $ translate of $C$, and fix a point $t \in C$.  Observe that
the intersection $C_t \cap W_1 \cap W_2$ is the point $w_1$.  We shall
follow the convention to indicate a rational function on $C_t$ by
using the same name given to the corresponding function on $C$.
  
Consider $K := W_1 \otimes f + W_2 \otimes f $.  $K$ is the 
{\it  basic hyperelliptic cycle} of \cite{collino-jag}, where it was
provedthat it is a non trivial indecomposable element of
$CH^{g}(J(C),1)_{\mathbb{Q}}$ for generic $C$. There it was shown that
the primitive contribution of the standard regulator image of $K$ does
not vanish by studying an infinitesimal invariant of Griffiths type
associated with the relevant normal function.  The following
proposition shows that $K$ is real regulator indecomposable. Here
$\tau_C$ is of the form considered in section \ref{subsec:zm}.

\begin{prop} 
  For $C$ a generic hyperelliptic
  curve of genus $g \geq 2$ there is $ \tau_C$ with $ [R(K) ,\tau_C]
  \ne 0 \, .$ \label{prop:nonzero}
\end{prop}
\noindent  By definition $[R(K) ,\tau] = 2 \int_{C} log|f|  \tau_{|C}$;
we shall prove that it is not trivial by means of a reduction process
to the case of elliptic curves.

%\subsubsection{}

Let $E_{\lambda}$ be the elliptic curve with affine equation $y^2 =
x(x-1)(x-{\lambda})$. Define $f_{\lambda} := x$ as a rational function
on $E_{\lambda}$. The form associated with the $\Theta$ divisor is
here $ \theta_{\lambda} = (i/2) {\omega_{\lambda} }
\wedge\bar{\omega}_{\lambda} $.  It is an invariant form on
$E_{\lambda}$ of volume one.  We write $I({\lambda}) :=
\int_{E_{\lambda}} \log|f_{\lambda}| \theta_{\lambda} $.
\begin{lem} 
  $I({\lambda})$  varies with ${\lambda}$.
\end{lem}
\begin{proof}
  Multiplication of $x$ by ${\lambda}^{-1}$ shows that $E_{\lambda}$
  and $E_{\lambda^{-1}}$ are isomorphic models of the same curve $E$.
  On $E$ the volume forms coincide, while $f_{\lambda} = {\lambda}
  f_{{\lambda}^{-1}}$. Thus
  \begin{equation*}
  I({\lambda})= \int_{E} \log|{\lambda}| \theta + I({\lambda}^{-1})
  = \log|{\lambda}|+ I({\lambda}^{-1}) \, ,
  \end{equation*}
\noindent hence $I({\lambda})$ cannot be constant.
\end{proof}

%\subsubsection{} 
We define $I(h,\tau) := \int_{C} \log|h| \tau_C $,
for any rational function $h$ on $C$.
\begin{lem} 
  If $C$ is a generic curve of genus $2$,
  then $I(f,\tau) \ne 0$. \label{lem:If}
\end{lem}
\begin{proof}
  We prove it for a bielliptic curve $C$ which is a double cover of
  $E_1:= E_{\lambda_1}$ and of $E_2:= E_{\lambda_2}$.  Consider the
  diagram
\begin{equation*}
\begin{CD}
  E_2 @< k_2 << C @> k_1 >> E_1
  \\
  @V f_2 VV @V f VV @V f_1 VV
  \\
  \mathbb P^1 @< h << \mathbb P^1 @> h >> \mathbb P^1
\end{CD} 
\end{equation*}
\noindent  
Here $f_i$ is ramified over $\{0,1,\infty, \lambda_i \}$, $h$ is the
double cover ramified over $ \lambda_1 $ and $ \lambda_2 $, and $f: C
\to \mathbb P^1$ is the hyperelliptic cover ramified at $
h^{-1}(\{0,1,\infty \}) $. On the range of $h$ we have already fixed a
standard parameter, we choose a standard parameter on the domain of
$h$ so that $0$ maps to $0$, and similarly for $1$ and for $\infty$.
In this manner $f$ is a well defined rational function on $C$, and we
denote by $\bar f$ its transform under the involution of $\mathbb P^1$
associated with $h$. Letting $g:=hf$, we see that $f \bar f
= cg$, $c$ a constant.

One can take $\tau_C$ to be the form $(k_1^{*}(\theta_{1}) -
k_2^{*}(\theta_{2}))$ on $J(C)$ and thus
\begin{equation*} 
  I(g,\tau_C) = \int_{C} \log|g| ( k_1^{*}(\theta_{1}) - 
  k_2^{*}(\theta_{2}) )  \, \neq  0 \, .
\end{equation*}
\noindent  
It then follows that $I(f, \tau_C) \neq 0 $ for the general bielliptic
curve C because
\begin{equation*}
  I(f, \tau_C) + I(\bar f, \tau_C) = I(g, \tau_C ) + log|c| \int_{C}
  \tau_C \, = I(g, \tau_C) \, . \qed 
\end{equation*}
\renewcommand{\qed}{} \end{proof}

\begin{proof}[Proof of Proposition \ref{prop:nonzero}]
  The proof for arbitrary genus is obtained by induction. Starting
  from a hyperelliptic curve $G \to \mathbb P^1$ and a double cover $
  h: \mathbb P^1 \to \mathbb P^1$ we construct the commutative diagram
\begin{equation*}
\begin{CD}
  C @> {\pi} >> G
  \\
  @V f_C VV @V f_G VV
  \\
  \mathbb P^1 @> h >> \mathbb P^1
\end{CD} 
\end{equation*}
Here $C$ is the normalization of the Cartesian product, hence $\pi$ is
branched at $4-2m$ points, where $m$ is the number of ramification
points of $h$ which coincide with points of ramification for $f_G$. We
have $g(C) = 2g(G) +1 -m$.

By induction, Proposition \ref{prop:nonzero} holds for $G$. It then
also holds for $C$ by using the arguments given for Lemma
\ref{lem:If}, where  we now take $\tau_C$ to be the form $\pi^{\ast}
\tau_G$, the lift to $J(C)$ of $\tau_G$.
\end{proof}

\subsection{Indecomposable elements on hyperelliptic Jacobians}
For a $g$ dimensional abelian  variety $A$,
the following decomposition of the
higher Chow groups is a consequence of the motivic decomposition
of the diagonal due to Beauville \cite{be2}, and
Deninger and Murre \cite{dm}:

\begin{equation} \label{eq:decomp}
  CH^k(A,m)_{\mathbb{Q}} = \bigoplus_s CH^k(A,m)_s
\end{equation}
Here $CH^k(A,m)_s$ is the subspace of
$CH^k(A,m)_{\mathbb{Q}}$ on which $[n]^*$ (resp. $[n]_*$) acts
by multiplication by $n^{2k-s}$ (resp. $n^{2g-2k+s}$).
The Fourier transform of Mukai and Beauville induces
isomorphisms:
\begin{equation} \label{eq:fourier}
  \mathcal{F}:CH^k(A,m)_s \stackrel{\cong}{\longrightarrow}
  CH^{g-k+s}(\hat{A},m)_s
\end{equation}
where $\hat{A}$ is the dual abelian variety.

For $\Theta$ a symmetric ample divisor, the motivic hard Lefschetz 
theorem of K\"unnemann  \cite{kunnemann-lefs}
implies that intersecting with powers of $\Theta$
gives isomorphisms:
\begin{equation}  
  \cdot \ \Theta ^{g+s-2k}: CH^{k}(\mathcal{A},m)_s \stackrel{\cong}{\longrightarrow}
  CH^{g+s-k}(\mathcal{A},m)_s \ , \  0 \leq 2k-s \leq g \, .  \label{eq:lefs}
\end{equation}
 
It follows from the definitions that 
the decomposition \eqref{eq:decomp} and the isomorphisms
\eqref{eq:fourier} and \eqref{eq:lefs} preserve decomposable cycles
as defined in \eqref{def:decomp}, hence
\begin{gather}
  CH^k_{\mathrm{ind}}(A,1)_{\mathbb{Q}} = \bigoplus_s
CH^k_{\mathrm{ind}}(A,1)_s 
  \label{eq:decomp2} \\
  \mathcal{F}: CH^k_{\mathrm{ind}}(A,1)_s \stackrel{\cong}{\longrightarrow}
  CH^{g-k+s}_{\mathrm{ind}}(A,1)_s \label{eq:fourier2} \\
  \cdot \ \Theta ^{g+s-2k}: CH^{k}_{\mathrm{ind}}(\mathcal{A},1)_s \stackrel{\cong}{\longrightarrow}
  CH^{g+s-k}_{\mathrm{ind}}(\mathcal{A},1)_s\, , \ \ 0 \leq 2k-s \leq g \, .  \label{eq:lefs2}
\end{gather}

\begin{prop}
  Let $A$ be a g-dimensional abelian variety. Then
  $CH^g_{\mathrm{ind}}(A,1)_s=0$ for $s<2$ or $s>g$. \label{prop:optimal}
\end{prop}
\begin{proof}
  We use equation \eqref{eq:fourier2} i.e.
  $\mathcal{F}(CH^g_{\mathrm{ind}}(A,1)_s) = CH^s_{\mathrm{ind}}(\hat{A},1)_s$.  If
  $s<1$, $CH^s(\hat{A},1)$ is itself zero. The action of $[n]^*$ on
  $CH^1(\hat{A},1) = \mathbf{C}^{*}$ is trivial hence
  $CH^1(\hat{A},1)_1$ is also zero.  We conclude the proof by
  observing that for any $g$-dimensional variety $X$, $CH^s(X,1)=0$
  for $s> g+1$ and $CH^{g+1}_{\mathrm{ind}}(X,1)=0$
\end{proof}

\begin{rem}
  A conjecture of C.~Voisin \cite{voisin-products} says that
  $CH^2_{\mathrm{ind}}(X,1)$ should be countable for any smooth projective
  variety $X$. For an abelian
  variety $A$, the injectivity of the rational regulator on
  $CH^k(A,1)_2$ would imply that $CH^k_{\mathrm{ind}}(A,1)_2$ is
  countable. If $g=2$, the  proposition shows that 
  then $CH^2_{\mathrm{ind}}(A,1)_{\mathbb{Q}}$ would also
  be countable.
\end{rem}

For the rest of this section, $C$ will be a hyperelliptic curve of genus $g$
and $K$
the basic hyperelliptic cycle in $CH^g(J(C),1)$.
\begin{lem}
  The component of $K$ in $CH^g_{\mathrm{ind}}(J(C),1)_s$
  is zero for all $s \neq 2$. Consequently, for any integer
  $n$, $[n]_*(K) = n^2K$ in $CH^g_{\mathrm{ind}}(J(C),1)_{\mathbb{Q}}$.
  \label{lemma:hyp}
\end{lem}
\begin{proof}
Consider the following copies of $C$
embedded in $C \times C$: $X_1 = C \times \{w_1\}$,  $X_2= \{w_1\} \times C$,
$X_3 =C \times \{w_2\}$,  $ X_4 =\{w_2\} \times C$, $X_5 =\Delta$,
and $X_6 = \Delta'$. Here $w_1$ and 
$w_2$ are two distinct Weierstrass points, $\Delta $ is the diagonal,
 and $\Delta'$
is the image of $C$ via the embedding $x \mapsto (x,\sigma(x))$,
with $\sigma$ the hyperelliptic involution. If $f$ is a
Weierstrass function with $div(f) = 2w_1 - 2w_2$ then one
easily checks that $Z = \sum_{i=1}^6 X_i \otimes f$ is
an element of $CH^2(C \times C, 1)$. If we use $w_1$ to
embed $C$ in $J(C)$, then the image of $Z$ in
$CH^g(J(C),1)$ is equal to $2K - (1/2)[2]_*(K)$.

The involution $\tau = (id,\sigma)$ of $C \times C$ preserves $Z$, hence
$Z$ must be the pullback of a cycle from 
$C \times C/\tau \cong C \times \mathbb{P}^1$.
$CH^2( C \times \mathbb{P}^1,1)$ is always decomposable for any curve $C$,
hence $Z$ must also be decomposable. Thus $4K = [2]_*(K)$
in $CH^g_{\mathrm{ind}}(J(C),1)_{\mathbb{Q}}$. This implies that all the components
of $K$ except the one in $CH^g_{\mathrm{ind}}(J(C),1)_2$ must be zero.
\end{proof}

For $T = (t_1, t_2, \ldots, t_m)$, $P=(p_1,p_2, \ldots, p_m)$ points
on $C^m$, 
let
$K_P(T) =  K*([t_1 - p_1] - e)*([t_2 - p_2] - e)*\ldots*([t_m - p_m] - e)$ $ \in
CH^g(J(C),1)$, where $e$ is the origin in $J(C)$. 
The following theorem is the main result of this section.
\begin{thm}
  If $C$ is a generic hyperelliptic curve
  and $T$ a generic point of $C^m$, $1 \leq m \leq g-2$, then
  $K_P(T)$ 
  is indecomposable and the set of all such cycles with $T$ varying and
  $P$ fixed forms an uncountable
  subset of $CH^g_{\mathrm{ind}}(J(C),1)_{\mathbb{Q}}$.
  Furthermore,
  $CH^k_{\mathrm{ind}}(J(C),1)_s$ is uncountable for $3 \leq s \leq k \leq g$.
  \label{thm:hypmain}
\end{thm}

\begin{proof}
  The first part of the theorem follows directly by combining
  Propositions \ref{prop:zm}, \ref{prop:hodgegroup} and \ref{prop:nonzero}.
  
  Let $I = \oplus_{s=1}^g CH^g(J(C))_s$ be the zero cycles of
  degree 0 (with $\mathbb{Q}$ coefficients).
  Then for $m \geq 1$,  $I^{*m} = \oplus_{s=m}^g CH^g(J(C))_s$
  and  if $T \in C^m$, then $K_P(T) 
  \in CH^g(J(C),1)_{\mathbb{Q}}*I^{*m}$.
  Using Proposition \ref{prop:optimal}, Lemma \ref{lemma:hyp}
  and the fact that the subscripts are
  additive under Pontryagin products,
  we see that if $m=g-2$ then the only nonzero component of $K_P(T)$
  in $CH^g_{\mathrm{ind}}(J(C),1)_{\mathbb{Q}}$ is the one in
  $CH^g_{\mathrm{ind}}(J(C),1)_g$.
  The first part of the theorem then implies that the image
  of $K*CH^g(J(C))_{g-2}$ in $CH^g_{\mathrm{ind}}(J(C),1)_g$
  is uncountable.
  We then use that $CH^g(J(C))_s =CH^g(J(C))_1^{*s}$, $s \geq 1$, and
  descending induction on $s$ starting from $s=g$ to show
  that the image of  $K*CH^g(J(C))_{s-2}$ in $CH^g_{\mathrm{ind}}(J(C),1)_{s}$
  is uncountable for $3 \leq s \leq g$.

  The statement for $k<g$ follows by using the Fourier transform
  and the motivic hard Lefschetz theorem:
  Letting $k=s$, we see by \eqref{eq:fourier2} that $CH^s_{\mathrm{ind}}(J(C),1)_s$
  is uncountable for $3 \leq s \leq g$. By \eqref{eq:lefs2},  intersection
  with $\Theta^{g-s}$ induces an isomorphism from 
  $CH^s_{\mathrm{ind}}(J(C),1)_s$ to  $CH^g_{\mathrm{ind}}(J(C),1)_s$,
  hence intersection with $\Theta^{k-s}$ must induce an injection
  from $CH^s_{\mathrm{ind}}(J(C),1)_s$ to $CH^k_{\mathrm{ind}}(J(C),1)_s$,
  $3 \leq s \leq k \leq g$. Since  $CH^s_{\mathrm{ind}}(J(C),1)_s$ is
  uncountable, it follows that  $CH^k_{\mathrm{ind}}(J(C),1)_s$ is
  also uncountable.
\end{proof}

\begin{rem}
  Proposition \ref{prop:optimal} shows that the theorem is optimal for
  $k=g$. It should also be optimal for $k<g$ --- the proof for
  $s>k$ still works but we do not know how to handle the $s<2$ case.
\end{rem}

\section{ Decomposability specializes}  \label{sec:decomp}
 
Our aim is to show that decomposability specializes. Consider a flat
and projective family $ \mathcal X \to A $, where $A$ is a smooth
curve, $ \mathcal X $ is non singular and so is $X_{ 0 }$, the fibre
over $ p_{0}$.
\begin{thm} If the restriction of an element  $ \mathcal Q \in
  CH^{d}( \mathcal X,1)_{\mathbb{Q}}$ to the generic fibre is
  decomposable then restriction of $ \mathcal Q $ to the central fibre
  is also decomposable. \label{thm:decomp}
\end{thm}

\subsection {} 
One can replace $ \mathcal {Q}$ by a multiple and 
$A$ by an open subset of a finite cover of the
original $A$ so that now the following holds:
\begin{assumption} Over  $U:= A - \{ p_{0} \} $ the restriction 
  $ \mathcal Q_{U} $ is equivalent in $ CH^{d}( \mathcal
  X_{U},1)$ to an element $\mathcal W := \sum \mathcal Z_j
  \otimes f_j $, where $\mathcal Z_j$ are irreducible subvarieties
  which intersect properly $X_{ 0 }$ and $ f_j $ are rational
  functions lifted from $A$ and regular on $U$.
\end{assumption}

We meet now the problem that $ \mathcal W$ may have a boundary $B$ on
$\mathcal X$, and then $B$ is supported on $X_{ 0 }$.  On the other
hand $ \mathcal Q $ is a cycle, and so it has no boundary. Since on $
\mathcal X_{U}$ the classes $ \mathcal Q $ and $ \mathcal W$ coincide
by hypothesis, then under the boundary map $ CH^{d}( \mathcal X_{U},1)
\to CH^{d-1}( X_ {0} ) $ the image of $ \mathcal W$ is trivial, that
is $B$ is rationally equivalent to 0 on the central fibre. Let $R$ be
a relation of rational equivalence on $X_{0}$ which kills $B$, then $
\mathcal W - R$ has no boundary and thus $\mathcal W - R -\mathcal Q $
is a cycle for $ CH^{d}( \mathcal X,1)$ whose class is represented by
a cycle $M$ coming from $CH^{d-1}( X_0,1) $.  We will show that $M$
and $R$ can be moved conveniently for our purposes, but to do this we
need to embed $ \mathcal X $ in a larger space $ \mathcal Y \to A$.
Our first result is that we can work on $ \mathcal Y$ and with a {\it
  cycle } $\mathcal W^{Y} := \sum \mathcal Z_j \otimes f_j $ as
described above. The improvement is that $\mathcal W^{Y}$ has no
boundary now and that it is equivalent to $ \mathcal Q $ on $ \mathcal
Y$.  The second issue is to compute the restriction of $\mathcal W^{Y}
$ to $X_ {0}$ so to check that it is indeed of decomposable type. The
difficulty that we meet is the fact that each $f_j $ may have a
boundary on the central fibre, we overcome this trouble by means of
some test curves.

\subsection {} We find useful to work for a while with
$G(T)$ the $K$-theory of coherent sheaves on a quasiprojective scheme
$T$ and to use the topological filtration $F_{m}¥G_{1}(T)$, this is
the subgroup generated by the images of $i_{Z,\ast }:G_{1}( Z ) \to
G_{1}(T)$, where $i_{Z}: Z \to T $ are closed subschemes of $T$ of
dimension at most $m$.  The relations with Bloch's
$CH^{p}(X,1)$ are now recalled.

We know from Soul\'e \cite{soule-operations}
that the Quillen coniveau spectral sequence
$$
E_{1}^{p,q} = {\coprod}_{x \in X_{p}} K_{-p-q}k(x) \Rightarrow
G_{-p-q}(X) $$
degenerates modulo torsion for $p+q \leq 2$, in
particular we need
$$
G_{1}(X)_{\mathbb Q} = \bigoplus_{p=0}^{d}E_{2}^{p,-p-1}\otimes
{\mathbb Q}\,$$
Moreover for $ m =0,1,2$ the coniveau filtration of the
spectral sequence coincides rationally with the $\gamma $ filtration.
Using Bloch's isomorphism
$$
\bigoplus_{i}gr^{i}_{\gamma} G_{m}(X)_{ \mathbb Q} \simeq
G_{m}(X)_{ \mathbb Q} \simeq \bigoplus_{i}CH^{i}(X,m)_{\mathbb{Q}}$$
one has then
$$
E_{2}^{p,-p-1}\otimes {\mathbb Q}\simeq CH^{p}(X,1)_{\mathbb{Q}}$$
\bigskip

Let $ i_{D}¥ : D \to T $ be the inclusion of an effective Cartier
divisor, then $ i^{ \ast }_{D}¥ G_{m}(T) \to G_{m}(D) $ is defined.
Next proposition, see 2-11 in \cite{quillen-ktheory1}, implies that if
each component of a subscheme $Z$ of dimension $m$ intersects properly
$D$, then $ i^{ \ast }_{D} ( i_{Z \ast }(G_{1}( Z ) )$ lands in
$F_{m-1}¥G_{1}(D)$.

\begin{prop}[Quillen] \label{prop:quillen}
Consider a cartesian diagram of quasiprojective schemes
$$\begin{CD}
   X' @> g' >> X
  \\
  @V f' VV @V f VV
  \\ 
  Y'  @> g >> Y
\end{CD} $$ 
Assume that $f$ 
is proper, that $g$ is of finite tor-dimension and that $Y'$ and $X$
are tor-independent over $Y$ (i.e. 
$ Tor _{i}^{ O_{Y,y } } (   O_{(Y',y')} ,  O_{(X,x) }  ) = 0 $ for  
$i    \ge  1 $)

Then 
$$ {g^{ \ast } } {f _{ \ast }} = f'_{ \ast }  {g'} ^{ \ast }  : G_{m}(X) \to G_{m}(Y') $$

\end{prop}

\subsection  {A larger space }
 
We come back to consider $   \mathcal X \to A $, by  gluing 
it  transversally  with $X \times A$ along the fibre over $p_{0}$ 
we construct $\mathcal  Y$, so that   $g: \mathcal  Y  \to  A$
has central fibre the scheme $Y$, whose reduced support is $X_{0}$.

The diagram 
$$\begin{CD}
   X  @>   >> \mathcal X   @.
  \\
  @V   VV @V   VV    @.
  \\ 
  Y   @>   >> \mathcal Y @. \, :=  \mathcal X \cup _{X \times \{ p_{0} \}} (X \times A)
\end{CD} $$
satisfies the hypotheses of Proposition \ref{prop:quillen}.
Indeed one has only to check vanishing of $Tor_{1}$,
because $Y$ is a Cartier divisor in $\mathcal Y$.
Let $t$ be a local parameter at $p_{0}$,
then  $t$ generates the ideal of $Y$,
and we see that   $Tor_{1} =0$ because 
$t$ is not a 0-divisor in the local rings of $\mathcal X$.

Applying devissage we find then:
$$\begin{CD}
   G_{q}(X)  @<   << G_{q}(\mathcal X)
  \\
  @V \simeq  VV @V   VV
  \\ 
 G_{q}(Y)   @<   << G_{q}( \mathcal Y)
\end{CD} $$
        We use next the setting above to check  that
the specialization map   $ G_{1}(\mathcal X) \to G_{1}(X) $
is compatible with the topological filtration $F_{m}$. The useful fact is that  
on  $\mathcal Y$  it is easy to move to general position with respect to $X$.
This way of moving turns out to be also the trick that we need to 
show that decomposability specializes. 

   We work by induction on $m$  and  consider   $ z $ in the graded quotient 
$F_{m+1}/F_{m} G_{1}(\mathcal X)$, then $ z $
is represented by a Quillen cycle $ \sum  Z_{i}\otimes f_{i} $, with  
$ \sum  div(f_{i}) = 0$ where $Z_{i}$ are irreducible subvarieties of 
dimension $m+1$. Now:

 \begin{lem}  On $ X \times A$  a Quillen  chain $Z \otimes f  $ is
equivalent to a chain whose support intersects properly  $ X \times \{ p_{0} \} $.
\end{lem}

\begin{proof}
The only problem is when  $Z  \subset  X \times \{ p_{0} \} $.
Consider the symbol 
$ \{ g, f \} $ on  $ Z   \times A $,
where $g $ is a rational function with simple 0 at  $p_0$
and therefore $div(g) =  p_0 + \sum m_{i}q_{i}$.
The  boundary of $(Z   \times A) \otimes  \{ g, f \} $ is 
 $  ( Z \times  \{ p_{0} \}) \otimes f + 
\sum   (Z \times  \{ q_{i}  \}) \otimes  f^{m_{i}}
- (div (f) \times A) \otimes g $,
and thus $( Z \times  \{ p_{0} \}) \otimes f$ is equivalent to a chain of the stated
kind.
\end{proof}

        In this way  the image of   $z$ in  $  F_{m+1}/F_{m}G_{1}(\mathcal Y)$
is supported on a scheme $T$ of dimension $m+1$ with the property that each component
of $T$
meets properly the central divisor $X$,
therefore the restriction of $z$  belongs to  $  F_{m}G_{1}(X) $
because it is supported on the intersection $T \cap X$.

        The same proof yields that the chain $M+R$
which was discussed above is equivalent   on $\mathcal Y$
to a sum of Quillen chains of the preceding type $(div (f) \times A) \otimes g $.
Notice that here we may have had to replace  $A$ by a smaller open set, with this proviso
we have then proved:
 \begin{prop}  The restriction of $Q$ to $X$ 
coincides with the restriction from $ \mathcal Y$ to  $Y$ of a cycle 
$\mathcal W^{Y} :=  \sum  \mathcal Z_j \otimes  f_j $,
where $\mathcal Z_j$ are irreducible varieties in $\mathcal Y $ and where $ f_j $ is pull back 
of a rational function of the same name from $A$ under the flat projection   $\mathcal Z_j \to A $.
The functions $ f_j $ are regular away from $p_0$.
\end{prop}
\subsection {Computing the restriction of $\mathcal W^{Y}$}

Let  $(\cup_{j} Z_j) \cap X_0 = \cup_{l}H_l$  be the  decomposition in irreducible components, then
we know from the preceding discussion that the restriction morphism $  CH^{d}( \mathcal Y,1)
\to  CH^{d}(X,1)$ sends  $\mathcal W^{Y}$
to a cycle supported on $\cup_{l}H_l$, and therefore to a cycle $ \sum _{l}  H_l \otimes  h_l$. 
Next proposition completes the proof of the theorem.
\begin{prop} The rational functions $ h_l$ are constant functions.
\end{prop}

\begin{proof}    
By taking general linear sections of $\mathcal Z = \cup \mathcal Z_j$
one may construct curves  $ B = \cup  B_j$ such that $ B \to A $ is finite.
Our program is to prove that at the points of intersection of $ B$
with $H_l$ the value of the relevant function $ h_l$ is constant,
independent of the chosen section $B$. 

We may assume that $B$ is smooth outside the 
inverse image of $p_0$. Let $Q$ be the element of $G_1(B)$ determined by 
$\mathcal W^{Y}$, and thus $Q$ is such that 
its restriction to the smooth part is given by functions $f_j$
on the components $ B_j$ of $B$, each $f_j$ being the pull-back of a rational
function with the same name on $A$. Using devissage  we restrict 
$Q$ to $G_1(p) \simeq \mathbb C^{*}$,
where $p$ is a point  in the inverse image of $p_0$. 
Then the claim is that for any such
point $p$ this restriction is equal to ($ \prod {f_j}^{n_j})(p_0)$.
Here the $n_j$'s may depend on $p$ and they are such
that the product doesn't have a pole or zero at $p_0$. Note that 
contravariant functoriality for local complete intersections 
yields here the equality  $ (\prod {f_j}^{n_j})(p_0)$ $=$
$ \prod (h_l(p) )$, where we set $h_l(p) = 1$ if $p \not \in H_l$.
This argument shows that at the points of intersection of curves like $B$
with $H_l$  the values of the functions $h_l$  
lie in a countable set, and therefore $h_l$ 
must be constant on the irreducible component $H_l$.

To prove the claim, we let $C$ be the normalization of an irreducible
component of the fibre product of all the $B_i$'s over $A$. $C$
is flat over $A$, so we can replace $C$ by $A$ and $B$ by the fibre
product of $A$ and $B$ over $C$. The advantage is that here we are
reduced to the case that there is only one point in the
inverse image of (the new) '$p$' because we now have sections,
and so the claim follows by functoriality from the simple case,
which we describe next.

The simple case is when  each
$B_i$ is an isomorphic copy of $A$ and moreover 
the scheme  $B:= \cup B_i$, is constructed by
gluing at one point $q$ (= $p_0$   on $A$).
Let the projection map be $ \pi : B \to  A$, then 
we write  $ U := A - \{p_0 \}$, and $ V_i = B_i - \{ q \}$,
and set $ V := \cup V_i = B -  \{ q \}$.
The map $ \pi : B \to  A$
is the identity on each component $B_i$.
Give rational functions $f_i$ on $B_i$,
($A$ is not necessarily complete, and we assume that the only zero or pole
of $f_i$ is at $q$)
and assume that $ \sum div (f_i) = 0 $.
Then $  \{f_i \}$ defines an element $  f $ say in 
$G_1(V)$, which comes from $G_1(B)$,
because it has boundary 0
in the exact sequence 
$G_{1}(B) \to G_{1}(V) \to G_{0}(q)$.
Consider now $p_0 $ as a  Cartier divisor on $A$ 
and let $ q^{*}$ be the scheme, which is the pull back of $p_0$ on  $B$.
The map $G_1(B)  \to G_1(q*) $ is well defined because
$ q^{*}$ is a local complete intersection in $B$. By devissage 
$G_1(q^{*}) = G_1(q) = \mathbb C ^{*} $;
the question is to understand  the value of the image of 
$ f $ in $G_1(q) = \mathbb C ^{*}$.
Each $f_i$ is a function $ f _{i A}$ on $A$. 
The answer is consider the product $ f _{A}  := \Pi  f _{i A}$,
then $ f _{A}$ is a rational function which is in fact regular
at $p_0$, because of our assumption, 
and then we have that the image value is $ f _{A}(p_0)$. 
The reason is once more the commutativity from
Proposition \ref{prop:quillen} (it is so below at $p_0$ and then it
must have been so at $q$).
\end{proof}

\section{The 4-configuration} \label{sec:4config}

For $C$ a generic hyperelliptic curve of genus $g \geq 3$, we have
constructed indecomposable elements in $CH^g(J(C),1)_{\mathbb{Q}}$ which
are in the kernel of the regulator map.  We shall now construct such
elements in $CH^g(J(C),1)_{\mathbb{Q}}$, for
more general curves, in particular for generic curves of
genus $3$ and $4$. The natural parameter space for our construction
will be a certain Hurwitz scheme of degree $4$ covers of 
$\mathbb{P}^1$. We shall show that one can specialize to a hyperelliptic
curve in such a way that our cycle specializes to $K- K_t$, allowing us to
use the results of the previous sections to prove indecomposability.

\subsection{Construction of the 4-configuration}

Let $C$ be a smooth projective curve of genus $g \geq 3$ and let $D =
(a_1 + a_2) -(b_1 + b_2)$ be a divisor on $C$ such that $[D] =
\epsilon$ is a point of order $2$ in $Pic(C)$. Let $f$ be a rational
function on $C$ such that $div(f) = 2D$. To this data we shall
associate a natural element of $CH^g(Pic^3(C),1)$ supported
on four copies of $C$.

Let $ i(y,z): C \to Pic^3(C) $ be the map $ i(y,z) (x) = x+y+z$,
$C(y,z) := i(y,z) (C)$, and let $ j : C \to Pic^3(C) $ be the map $ j
(x) = -x+2(a_1+a_2)$, $G := j(C)$.  We consider $ C(a_1,a_2)$,
$C(b_1,b_2)$ and $G$. Translation by ${\epsilon}$ maps $ C(a_1,a_2)$
to $C(b_1,b_1)$, and we let $G_{\epsilon}$  be the image of $G$. We
shall use the convention that $f$ represents the rational function on
each of the preceding curves which maps to $f$ under the chosen
isomorphism with $C$. Thus, we set $ Z_1 := C(a_1,a_2) \otimes f$, $
Z_2 := G \otimes f $, $ Z_3 := C(b_1,b_2) \otimes f$, $ Z_4 :=
G_{\epsilon} \otimes f$.

\begin{prop} 
  $Z_D := \sum_{i=1}^{4} (-1)^i Z_i$ is a higher cycle.
\end{prop}

Note that $Z_D$ also depends on the choice of the rational function
$f$, but we can neglect this since multiplying $f$ by an element of
${\mathbb C}^{*}$ amounts to the addition of a decomposable
element to $Z_D$.

\begin{proof}
  The only possible difficulty is to see where
  the curves intersect. $ C(a_1,a_2)$ intersects $G $ in two points;
  on both curves the points come from $a_1$ and $a_2$ under the
  isomorphism with $C$, but the point which comes from $a_1$ in $G$
  comes from $a_2$ in $ C(a_1,a_2)$ and conversely. A similar
  statement holds for $C(b_1,b_2) \cap G$, and then intersections with
  $G_{\epsilon}$ can be recovered by using ${\epsilon}$-symmetry.
  Note that if $C$ is not hyperelliptic then $C(a_1,a_2)\cap
  C(b_1,b_2) = \emptyset $ and $ G \cap G_{\epsilon} = \emptyset $.
\end{proof}
\begin{rem}
By translation of $Z_D$ by a a zero cycle $\xi$ of degree $-3$ on $C$,
we obtain
an element $Z_{D,\xi}$ of  $CH^g(J(C),1)$.
The component of $Z_{D,\xi}$ in $CH^g(J(C),1)_2$ is always zero,
hence it is in the kernel of the (rational) regulator map.
To see this, first note that this component is preserved by
translation since $CH^g(J(C),1)_s = 0$ for $s < 2$ (c.f. proof of
Proposition \ref{prop:optimal}), hence is independent of $\xi$.
One easily checks that $[-1]^*(Z_{D,-3a_1}) = - Z_{D,-3a_2}$,
This implies the desired fact, since the action of $[-1]^*$
on  $CH^g(J(C),1)_2$ is trivial.
\end{rem}

\subsection{A Hurwitz scheme}
In order to construct the specialization needed to prove
indecomposability, we shall use certain Hurwitz schemes parametrizing
the ramified coverings of $\mathbb{P}^1$ corresponding to the
functions $f$ as above. We refer the reader to 
\cite[Section 1]{fulton-hurwitz} for the basic facts about Hurwitz
schemes.

Let $g \geq 2$ be an even integer and $H_g$ be the Hurwitz scheme
whose points correspond to degree 4 covers of $\mathbb{P}^1$ branched
over $n = 2g + 4$ distinct points, such that the inverse image of
$n-2$ of these points consists of three  points and the inverse image of
each of the remaining points consists of two points, each of
ramification degree $2$.  From the Riemann-Hurwitz formula it follows
that such a cover $C \to \mathbb{P}^1$ is of genus $g$.

\begin{prop}
  For each $g \geq 2 $, $H_g$ consists of two components.
  \label{prop:twocomp}
\end{prop}

\begin{proof}
  Let $n = 2g + 4$.  Recall that we have a finite etale map $\delta:
  H_g \to \mathbb{P}^n - \Delta$, where $\mathbb{P}^n$ is thought of
  as ${(\mathbb{P}^1)}^{(n)}$ and $\Delta$ is the discriminant locus.
  Let $P = \{p_1,p_2, \ldots, p_n\}$ be an element of $\mathbb{P}^n -
  \Delta$ and let $x \in \mathbb{P}^1$ be distinct from the $p_i$'s.
  We may choose loops $\sigma_i$ based at $x$ and going around $p_i$
  in such a way so that $\sigma_1,\sigma_2, \ldots \sigma_n $ generate
  $G = \pi_1(\mathbb{P}^1 \backslash \{p_1,p_2,\ldots, p_n\}, x)$ with
  the single relation as $\sigma_1\sigma_2 \cdots \sigma_n = 1$.  Now
  degree $4$ covers of $\mathbb{P}^1$ branched over $P$ correspond to
  transitive representations of $G$ in the symmetric group $\Sigma_4$.
  For the covers to correspond to points of $H_g$, the images of $n-2$
  of the $\sigma_i$'s must be transpositions and the images of the
  other two must be products of $2$ disjoint transpositions. Two such
  representations give isomorphic covers if they differ by an inner
  automorphism of $\Sigma_4$. Thus $\delta^{-1}(P)$ can be identified
  with classes of $n$-tuples $(s_1,s_2, \ldots, s_n)$ of elements of
  $\Sigma_4$ , all but two of the $s_i$'s being transpositions, the
  remaining two being products of two disjoint transpositions, and
  $s_1 s_2 \cdots s_n =1$. Two such $n$-tuples are identified if they
  differ by coordinatewise conjugation by an element of $\Sigma_4$.
  
  The action of the monodromy on $\delta^{-1}(P)$ contains elements
  $\Gamma_i$, $1 \leq i \leq n-1$, which act on $n$-tuples as above
  by:
  \[ \Gamma_i(s_1,s_2, \ldots, s_n) =
  (s_1, \ldots, s_{i-1}, s_i s_{i+1}s_i^{-1}, s_i, \ldots, s_n)\, .
  \]
  To show that $H_g$ has two components we use the $\Gamma_i$'s to
  prove that the monodromy action on $\delta^{-1}(P)$ has two orbits.
  The proof is a case by case analysis; we shall describe the main
  steps and leave some simple verifications to the reader. 
  Let
  $t_{i,j}$, $ i \neq j$, denote the transposition which switches $i$
  and $j$ and let $v_1$, $v_2$ and $v_3$ be the permutations $(1 \ 
  2)(3 \ 4)$, $(1 \ 3)(2 \ 4)$ and $(1 \ 4)(2 \ 3)$ respectively. Let
  $V = \{v_1,v_2,v_3\}$.
  Using the action of the $\Gamma_i$'s one sees that each orbit
  contains an $n$-tuple $S$ such that $s_{n-1}$ and $s_{n}$ are in
  $V$.  Upto conjugation, we may assume that $s_n = v_1$ and $s_{n-1}
  = v_1$ or $v_2$.  Suppose $s_{n-1} = v_2$, hence $s_1s_2 \cdots
  s_{n-2} = v_3$.  Let $K$ be the subgroup of $\Sigma_4$ generated by
  $s_1,s_2, \ldots s_{n-2}$.  Using the fact that $v_3 \in K$ and $K$
  is generated by transpositions we analyze the two possibilities for
  the action of $K$ on $\{1,2,3,4\}$ (i) transitive and (ii)
  intransitive.  In both cases ((i) requires some computations with
  the $\Gamma_i$'s) we conclude that we may assume that at least one
  of $t_{1,4}$ or $t_{2,3}$ occurs among the $s_i$'s (without changing
  $s_{n-1}$ and $s_n$).  Again using the action, we may assume that $i
  = n-2$.  Then replacing $S$, by $\Gamma_{n-2}(\Gamma_{n-1}(S))$ we
  may also assume that $s_{n-1}= v_1$.
  
  As above, let $K$ be the subgroup of $\Sigma_4$ generated by
  $s_1,s_2, \ldots s_{n-2}$ and note that now $s_1s_2 \cdots s_{n-2} =
  1$.  We now consider four possibilities for the action of $K$ on $
  \{1,2,3,4\}$.
\begin{enumerate}
\item The action is transitive.
\item The action has a unique fixed point.
\item The action has two fixed points.
\item The action has no fixed points but fixes two disjoint subsets of
  two elements each.
\end{enumerate}

\noindent (i) Since $s_1s_2 \cdots s_{n-2} = 1$, it follows by a
result of Clebsch (see \cite[p.547]{fulton-hurwitz}) that we may
assume
\[
(s_1,s_2, \ldots s_{n-2}) =(t_{1,2},t_{1,2},t_{1,3},t_{1,3},t_{1,4},
\ldots t_{1,4})\, .
\]
Note that a priori we may also have to use an inner automorphism of
$\Sigma_4$ in order to achieve this, but it is easy to check using the
action of the $\Gamma_i$'s that we may choose an automorphism
preserving $v_1$.  Hence we may assume that
\[
S=(t_{1,2},t_{1,2},t_{1,3},t_{1,3},t_{1,4}, \ldots t_{1,4},v_1,v_1)\, .
\]
\noindent (ii) Without loss of generality, we may assume that the fixed point
is $4$. Again, by Clebsch' result we may assume that
\[
S=(t_{1,2},t_{1,2},t_{1,3}, \ldots t_{1,3},v_1,v_1)\, .
\]
Moreover, it is easy to check that if we replace $v_1$ by $v_2$ or
$v_3$, we stay in the same orbit, hence all elements of type (ii) lie
in one orbit. Further, note that using the $\Gamma_i$'s , we can
switch pairs of adjacent transpositions, i.e.
\[
(t_{1,2},t_{1,2},t_{1,3},t_{1,3}) \sim
(t_{1,3},t_{1,3},t_{1,2},t_{1,2})\, ,
\]
so we may replace $S$ by $(t_{1,3}, \ldots,
t_{1,3},t_{1,2},t_{1,2},t_{1,3},t_{1,3},v_1,v_1)$.  One checks that
\begin{multline*}
  \Gamma_{n-2}(\Gamma_{n-3}(\Gamma_{n-4}(\Gamma_{n-5}(\Gamma_{n-4}(\Gamma_{n-5}(\Gamma_{n-4}(\Gamma_{n-3}(\Gamma_{n-2}(\Gamma_{n-4} (S)))))))))) \\
  = (t_{1,3}, \ldots, t_{1,3},t_{2,4},t_{1,2},t_{1,2},t_{2,4},v_1,v_1)\, .
\end{multline*}
Since we have assumed that $n \geq 8$, it follows that the $K$
corresponding to this element acts transitively on $ \{1,2,3,4\}$.
Hence elements of types (i) and (ii) lie in the same orbit. \\
%\vspace{3pt}
\noindent (iii) In this case  we must have
$S=(t,t,\ldots,t,v_1,v_1)$ where $t$ is one of $t_{1,3}$, $t_{2,3}$,
$t_{1,4}$, $t_{2,4}$, since the representation is assumed to be
transitive.  It is clear that all choices are conjugate by an element
of $\Sigma_4$ which fixes $v_1$, hence all elements of type (iii) are
in the same orbit. 

Assume, without loss of generality that $t=t_{1,3}$.  Then
\[\Gamma_{n-2}(\Gamma_{n-3}(\Gamma_{n-3}(\Gamma_{n-2}(S))))
= (t_{1,3},t_{1,3}, \ldots ,t_{1,3},t_{2,4},t_{2,4},v_1,v_1)\, .
\]
Since $t_{1,3}$ and $t_{2,4}$ commute, $\Gamma_i$, for $1 \leq i \leq
n-2$, just switches $s_i$ and $s_{i+1}$. Thus we may repeat the above
procedure and deduce that $S$ is in the same orbit as elements of the
form $(s_1, s_2, \ldots, s_{n-2}, v_1,v_2)$ where each $s_i$ for $1
\leq i \leq n-2$, is either $t_{1,3}$ or $t_{2,4}$, there being an
even number of both kinds.\\
%\vspace{3pt}
\noindent (iv) Elements of type (iv) are precisely those considered in
the previous paragraph, so all such elements lie in the same orbit as
elements of type (iii).

We thus see that $H_g$ has at most two components.  Now observe that
covers of type (iii) have an automorphism of order $2$ commuting with
the covering map, corresponding to the representation of $G$ in
$\Sigma_2$ induced from the original representation in $\Sigma_4$ by
identifying $1$ with $2$ and $3$ with $4$.  Covers of type (i) have no
automorphisms commuting with the covering map, hence $H_g$ has
precisely two components.
\end{proof}

\subsection{Indecomposability of the generic $4$-configuration}

Let $H_g'$ be the component of $H_g$ corresponding to covers without
automorphisms. Then there exists a universal family of curves $\psi:
\mathcal{C}_g \to H_g'$ and a corresponding universal degree $4$ map
$\pi : \mathcal{C}_g \to H_g' \times \mathbb{P}^1$.

We continue using the same notation as in the proof of Proposition
\ref{prop:twocomp}, so  $ P = \{p_1$, $p_2$, \ldots, $p_n\}$, 
$p_i \in \mathbb{P}^1$,
and let $p$ be any other point of $\mathbb{P}^1$.  Let $X$ be
the curve in $\mathbb{P}^n = (\mathbb{P})^{(n)}$
 with points $\{(1-t)p_1 + tp, (1-t)p_2 +
tp, p_3, p_4, \ldots p_n\}$, $t \in \mathbb{C}$. By the previous
proposition, there exists an element of $\delta^{-1}(P)$ such that the
monodromy representation of the corresponding cover is given by the
$n$-tuple
%$S= ((2 3), (2 3),(1 2),(1 2), \ldots, (1 2), (1 2) (3 4), (1 2) (3 4))$ 
$S = (t_{2,3},t_{2,3}, t_{1,2}, t_{1,2}, \ldots, t_{1,2}, v_1, v_1)$
of elements of $\Sigma_4$.  Let $Y'$ be the component of
$\delta^{-1}(X)$ containing this point and let $Y$ be the
normalization of $X$ in the function field of $Y'$.  Then $\mathcal{C}
= \psi^{-1}(Y')$ maps to $X \times \mathbb{P}^1$ by the map $\pi$ and we
let $\overline{\mathcal{C}}$ be the normalization of $X \times
\mathbb{P}^1$ in the function field of $\mathcal{C}$. Clearly
$\psi|_{\mathcal{C}}$ induces a flat and projective morphism
$\overline{\psi}:\overline{\mathcal{C}} \to Y$. Let $y_0$ be a point of $Y$ lying
above the point $\{p, p, p_3, p_4, \ldots p_n\}$ of $X$.

\begin{lem}
  The fibre of the map $\overline{\psi}: \overline{\mathcal{C}} \to Y$ over $y_0$
  has two smooth components, $C_1$ and $C_2$, meeting transversally in
  a single point lying over $p$. The map from $C_1$ to $\mathbb{P}^1$
  has degree $2$ and is branched over $p_3$, $p_4$, \ldots, $p_n$,
  while the map from $C_2$ to $\mathbb{P}^1$ is branched over
  $p_{n-1}$ and $p_n$. In particular, $C_1$ is hyperelliptic of genus $
  g = (n-4)/2$ while $C_2$ is of genus $0$. \label{lemma:badfibre}
\end{lem}

\begin{proof}
  Let $B$ be small topological disc in $Y$ containing $y_0$ and such
  that the only point where the map from $B$ to $X$ is ramified is
  $y_0$. By construction of $X$, if $B$ is small enough then the
  ramification locus of the map $\pi$ intersected with $(B - \{y_0\})
  \times \mathbb{P}^1 $ consists of $n$ disjoint punctured discs, each
  mapping isomorphically to $B - \{y_0\}$.  The closures of all these
  discs remain disjoint in $B \times \mathbb{P}^1 $ except for two of
  them which meet at $(y_0, a)$. By construction, the product of the
  local monodromies around these two discs is trivial, hence the map
  $\pi|_{\overline{\psi}^{-1}(B)}: \overline{\psi}^{-1}(B) \to B \times \mathbb{P}^1$
  is unramified over the complement of the closures of all the above
  discs.
  
  To complete the proof we examine the induced monodromy
  representation of $\pi_1(\mathbb{P}^1 \backslash \{a,
  p_3,p_4,\ldots, p_n\}, x)$ in $\Sigma_4$. The local monodromy around
  $a$ is equal to the product of the local monodromies around $p_1$
  and $p_2$ of the original representation used in the construction of
  $Y$ and the local monodromies around the other points are the same
  as those in the original representation. Then we see that the
  representation is no longer transitive but breaks up into two
  representations, each of degree $2$.  The conclusions of the lemma
  follow immediately on inspection of these two representations.
\end{proof}

For any point $h$ of $H_g'$, the covering map $C = (\mathcal{C}_g)_h
\to \mathbb{P}^1$ gives us points $a_1$, $a_2$, $b_1$, $b_2$ on $C$
such that $2[a_1 + a_2] = 2[b_1 + b_2]$. Choosing a function $f$ on
$C$ such that $div(f) = 2(a_1 + a_2) - 2(b_1 + b_2)$ allows us to
define a $4$-configuration $Z$ in $Pic^3(C)$.  Note that $Z$ is well
defined upto sign in $CH^g_{\mathrm{ind}}(Pic^3(C), 1)$.

\begin{thm}
  The $4$-configuration corresponding to a generic point of $H_g'$ is
  indecomposable for all $g \geq 3$.  Moreover $CH^k_{\mathrm{ind}}(J(C),1)_s$
  is uncountable for $4 \leq s \leq k \leq g$.
\label{thm:indec4}
\end{thm}

\begin{proof}
  Let $\overline{\mathcal{C}} \to Y$ be the family constructed in the
  discussion preceding Lemma \ref{lemma:badfibre}.  We blow down the
  genus $0$ curve $C_2$ in the fibre over $y_0$ and call the resulting
  family of curves $\mathcal{C}'$. By replacing $Y$ by a Zariski open
  subset we may assume that all fibres are smooth and then by
  replacing $Y$ by a finite cover we may also assume that branch locus
  of the map $\overline{\mathcal{C}} \to Y \times \mathbb{P}^1$ is a union
  of sections.  Finally, by replacing $Y$ by a further open subset
  (containing $y_0$) we may assume that there exist a rational
  function $F$ on $\mathcal{C}'$ such that $div(F) = 2(A_1 + A_2) -
  2(B_1 + B_2) = 2\mathcal{D}$, where $\mathcal{D}$ restricts to the
  divisor used to define the $4$-configuration on each fibre.
  
  $F$ allows us to construct an element $\mathcal{Z}$ of
  $CH^g(Pic^3(\mathcal{C}'/Y), 1)$ which restricts to the
  $4$-configuration on each fibre.  By Lemma \ref{lemma:badfibre}, we
  see that upto labelling the restrictions of $A_1$, $B_1$, $A_2$,
  $B_2$ to $C_1 = {\mathcal{C}'}_{y_0}$ must be $w_1$, $w_2$, $t$,
  $t$, where $w_1$ and $w_2$ are distinct Weierstrass points on $C_1$
  and $t$ is a point lying over the point $a \in \mathbb{P}^1$.  It is
  then clear that $f_{y_0} = F|_{C_1}$ has divisor $2w_1 - 2w_1$ and
  is therefore a Weierstrass function. One easily checks that the
  element of $CH^g(J(C_1), 1)$ obtained by translating
  $\mathcal{Z}_{y_0}$ by $-3w_1$ is equal to $K - K_t$, where $K$ is
  the basic hyperelliptic cycle and the subscript denotes translation.
  Since $p_1$, $p_2$, \ldots, $p_n$ and $p$ are arbitrary points in
  $\mathbb{P}^1$, it follows that we may assume that $C_1$ is a
  generic hyperelliptic curve of genus $g$ and $t$  a generic point
  on $C_1$. By Theorem \ref{thm:hypmain} it follows that
  $\mathcal{Z}_{y_0}$ is indecomposable and then by Theorem
  \ref{thm:decomp} it follows that $\mathcal{Z}$ restricted to a
  generic fibre is also indecomposable.  

  The second statement follows
  by considering the Pontryagin product of the
  $4$-configuration with the zero cycles $(t_1 - t_0)*(t_2 - t_0)*
  \cdots *(t_r - t_0)$, $0 \leq r \leq g -3$,
  where  $(t_0,t_1, \ldots, t_r)$
  is a generic point of $C^{r+1}$.
  Taking $r = g-3 > 0$, we see from the proof of Theorem \ref{thm:hypmain}
  that the components of elements of $Z*I^{*(g-3)}$
  give uncountably many elements
  of $CH^g_{\mathrm{ind}}(J(C),1)_g$ (note that specialization
  preserves the decompositions). We then deduce that 
  $CH^g_{\mathrm{ind}}(J(C),1)_s$ is uncountable for $4 \leq s \leq g$,
  as well as the  statement for $k < g$, in the same way as
  in Theorem \ref{thm:hypmain}.
\end{proof}

Consider the natural map $\tau_g:H_g' \to \mathcal{M}_g$.  Since the
map $C^{(2)} \times C^{(2)} \to J(C)$ given by $(\{a_1,a_2\},
\{b_1,b_2\}) \mapsto [a_1 + a_2 - (b_1 + b_2)]$ is surjective for any
curve of genus $g \leq 4$, by counting parameters we see that $\tau_g$
is dominant for $g = 3, 4$.  For $g > 4$ 
a degeneration argument shows that the image has dimension $2g + 1$.
Note that the (non-empty) fibres of $\tau_g$ are at least $3$ 
$=dim(Aut(\mathbb{P}^1))$ dimensional and the generic fibre
of $\tau_3$ is $4$ dimensional.

\begin{cor}
  For $C$ a generic curve of genus $g= 3,4$, $CH^g_{\mathrm{ind}}(J(C),1)_g$ is
  uncountable.
\end{cor}

\begin{proof}
  The statement for $g=4$ follows  from the theorem.  For
  $g=3$, the generic curve has a $1$-parameter family of
  $4$-configurations. Consider a family of genus 3 curves with special
  fibre a hyperelliptic curve as in the theorem. By varying the
  4-configuration in the Jacobian of the generic fibre in the
  $1$-parameter family, we obtain as specializations cycles of the
  form $K - K_t$ where now $t$ varies in a $1$-parameter family. By
  Theorem \ref{thm:hypmain}, it follows that the set of
  specializations is an uncountable subset of
  $CH^3_{\mathrm{ind}}(J(C_1),1)_3$. Hence by Theorem \ref{thm:decomp} the
  $4$-configurations also form an uncountable subset of
  $CH^3_{\mathrm{ind}}(J(C),1)_3$.
\end{proof}

\section{Higher Chow cycles on self products of a curve} \label{sec:products}

\subsection{Construction of the cycles}
The indecomposability of the basic hyperelliptic cycle and the
4-configuration suggest that to construct indecomposable elements in
$CH^g(J(C),1)$ for more general Jacobians, we should consider more
general divisors $D$ on $C$, with $[D]$ of order 2 in $Pic(C)$. Let
$D$ be such a divisor, so $D = \sum_{i=1}^n a_i -\sum_{i=1}^n b_i$ for
some $n > 0$ and distinct points $a_i, b_i \in C$, and there exists a
rational function $f$ on $C$ with $div(f) = 2D$.  Note that if $ 0 < g
\leq 2n$ then there always exist such divisors on a general curve, and
if $g \geq 2n$ then there exists a $2g + 2n - 3$ dimensional family of
curves of genus $g$ which have such divisors.

Let $E$ be the set of embeddings of $C$ in $C^{n+1}$ such that the
composition of any $e \in E$ with the projection to any factor is
either the identity map or is constant with image equal to one of the
$a_i$'s or $b_i$'s.  We have the following group actions on the set
$E$.
\begin{enumerate}
\item $\Sigma_{n+1}$ acts by permuting the factors.
\item $\Sigma{}_n$ acts by permuting the $a_i$'s.
\item $\Sigma{}_n$ acts by permuting the $b_i$'s.
\item $\mathbb{Z}/2$ acts by switching $a_i$ and $b_i$ for all $i$.
\end{enumerate}
Let $f_e$ be the function $f$ considered as a function on $e(C)$.
\begin{prop} 
  The set of higher Chow $1$-cycles of the form $\sum_{e \in E} m_e
  (e(C) \otimes f_e)$, $e \in E$, $m_e \in \mathbb{Q}$ invariant under
  all the above group actions on $E$, is a positive dimensional vector
  space for all $n >0$. \label{prop:products}
\end{prop}

\begin{proof}%[Sketch of proof]
  By a partition of a positive integer $n$ we shall mean a tuple of
  non-increasing positive integers $\alpha= (n_1,n_2, \ldots, n_r)$
  such that $|\alpha| := \sum_{i=1}^r n_i = n$ and by a partition of
  zero we mean the empty tuple $()$. If $\alpha$ is as above and $m$
  is a positive integer, we let $\alpha + m$ be the partition of $n +
  m$ obtained by reordering $(n_1,n_2, \ldots, n_r,m)$.
  
  The orbits of the action generated by all the above groups on $E$
  are in $1$-$1$ correspondence with the set
\[
\mathcal{E} = \big\{ \{\alpha,\beta\} \ | \ \alpha \text{ a partition
  of } i, \ \beta \text{ a partition of } j,\ i + j \leq n\big\} \, .
\]
Here the elements of $\mathcal{E}$ (and $\mathcal{P}$ below) are
viewed as unordered pairs.

The boundary of an element $\sum_{e \in E} m_e (e(C) \otimes f_e)$ is
supported on points of $C^{n+1}$, each of whose coordinates is one of
the $a_i$'s or $b_i$'s. The orbits of the set of all such points under
the group action are in $1$-$1$ correspondence with the set
\[
\mathcal{P} = \big\{ \{\alpha,\beta\}| \alpha \text{ a partition of }
i, \ \beta \text{ a partition of } j, i + j = n + 1\} \backslash
\{\{(1, \ldots,1), ()\}\big\} \, .
\]
$\{(1, \ldots,1), ()\}$ is not included because the number of $a_i$'s
and $b_i$'s is $n$.

Let $\mathcal{R} \subset \mathcal{E} \times \mathcal{P}$ be the set
\[
\bigcup_{\{\alpha,\beta\} \in \mathcal{E} } \! \!
\big\{(\{\alpha,\beta\},\{\alpha,\beta + (n + 1 -|\alpha| -
|\beta|)\}), (\{\alpha,\beta\},\{\alpha+ (n + 1 -|\alpha| -
|\beta|),\beta \})\big\} \, .
\]
Consider the projection $p_2:\mathcal{R} \to \mathcal{P}$ and observe
that $|p_2^{-1}(\{\alpha,\beta\})| \geq 2$ for all $\{\alpha,\beta\}
\in \mathcal{P}$ except for those of the form (i) $\alpha = (d,d,
\ldots, d)$, $\beta=()$, with $d > 1$ a divisor of $n + 1$ and (ii)
$\alpha = \beta = (d,d, \ldots, d)$ with $n$ odd and $d$ a divisor of
$(n+1)/2$.  For both these cases $|p_2^{-1}(\{\alpha,\beta\})| = 1$.
Using the fact that for any integers $d, m$, if $m/2 < d < m$ then $d$
cannot divide $m$, one checks that if $n>3$, then
$\big|\{\{\alpha,\beta\} \in \mathcal{P}\ | \ 
|p_2^{-1}(\{\alpha,\beta\})| = 1\}\big| \leq n$.

On the other hand, $|p_1^{-1}(\{\alpha,\beta\})| = 2$ for all
$\{\alpha,\beta\} \in \mathcal{E}$ except for those with $\alpha =
\beta$ or $\alpha = (1,1, \ldots,1)$, $|\alpha| = n$, and $\beta =
()$.  For both these cases $|p_1^{-1}(\{\alpha,\beta\})| = 1$.  If
$\Pi(m)$ denotes the number of partitions of $m$, one sees that the
number of such elements is $1 + \sum_{i=0}^{[n/2]} \Pi(i)$ which is
greater than $n$ for all $n > 5$.

Elements of $\mathcal{P}$ give us sufficient relations among the
$m_e$'s for the boundary of $\sum_{e \in E} m_e (e(C) \otimes f_e)$ to
be zero, but if $n$ is odd then elements of the form
$\{\alpha,\alpha\}$, $|\alpha| = (n+1)/2$, give trivial relations
because of symmetry.  Using this observation along with explicit
computations for $1 \leq n \leq 5$, one sees that the number of
relations is always strictly less than the number of variables i.e.
$|\mathcal{E}|$, hence the space of invariant cycles is positive
dimensional.
\end{proof}

Let $V(C,f)$ be the space of all invariant higher Chow 1-cycles as in
the proposition. By abuse of notation we shall also denote this space by
$V(C,D)$.  For $n=1$ and $2$, $V(C,D)$ is of rank $1$, but it is of
rank $>1$ for all $n>2$.

\subsection{Indecomposability of the product cycles via specialization}
Proposition \ref{prop:products}
provides us with  elements of $CH^{n+1}(C^{n+1}, 1)_{\mathbb{Q}}$.
After pushforward by the natural map $C^{n+1} \to Pic^{n+1}(C)$ and
translation by a divisor of degree $-(n+1)$, we obtain elements of
$CH^g(J(C),1)_{\mathbb{Q}}$.  We now outline a method which should
allow one to prove indecomposability of the image of a general element
of $V(C,D)$ in $CH^g(J(C),1)_{\mathbb{Q}}$, for $C$ generic among
curves of genus $g \geq n + 1 > 2$, having such divisors $D$.

Let $\mathcal{M}_g$ be the moduli space of genus $g$ curves 
(pointed curves if $g=1$) with level
$2m$ structure for some $m \geq 3$ and let $\mathcal{C}_g$ be the
universal family of curves over $\mathcal{M}_g$.  Consider the map
$\sigma:Sym^n(\mathcal{C}_g/\mathcal{M}_g) \times
Sym^n(\mathcal{C}_g/\mathcal{M}_g) \to J(\mathcal{C}_g/\mathcal{M}_g)$
given on the fibres over $\mathcal{M}_g$ by
$(\{a_1,a_2,\ldots,a_n\},\{b_1,b_2,\ldots,b_n\}) \mapsto [\sum_n a_i -
\sum_n b_i]$ and let $X_{g,n} = \sigma^{-1}(S_2)$, where $S_2$ is the
union of sections of $J(\mathcal{C}_g/\mathcal{M}_g)$ corresponding to
points of order $2$ on each fibre. Let $\tau$ be the composite
of $\sigma$ with the natural map from $J(\mathcal{C}_g/\mathcal{M}_g)$
to $\mathcal{M}_g$.

\begin{lem}
  Let $n> 1$ and $g > 0$. Then there exists an irreducible component
  $Y_{g,n}$ of $X_{g,n}$ with the following properties:
  \begin{enumerate}
  \item $dim(Y_{g,n}) = 2g + 2n -3$.
  \item If $g \leq 2n$, then the map $\tau|_{Y_{g,n}}:Y_{g,n} \to
    \mathcal{M}_g$ is dominant and if $g \geq 2n$, then
    $dim(\tau(Y_{g,n})) = 2g + 2n -3$.
  \item In the fibre of $\tau|_{Y_{g,n}}$ over a point on
    $\mathcal{M}_g$ corresponding to a generic hyperelliptic curve
    $C'$, there exist points $ (\{w_1,t_1,t_2, \ldots,
    t_{n-1}\},\{w_2,t_1,t_2, \ldots, t_{n-1}\})$, where $w_1$ and
    $w_2$ are distinct Weierstrass points on $C'$, and $(t_1,t_2,
    \ldots,t_{n-1})$ is a generic point of ${C'}^{n-1}$.
  \end{enumerate} \label{lemma:spec}
\end{lem}

\begin{proof}[Outline of proof]
  The lemma can be proved by considering a suitable Hurwitz scheme as
  in section \ref{sec:4config}. We do not know the number of
  components  if $n>2$, but a suitable choice of
  the monodromy representation allows us to single out a component
  which gives rise to the desired specializations
  (c.f. the discussion before Lemma \ref{lemma:badfibre}).
  For example, if
  $n=3$ we would consider the representation in $\Sigma_6$
  corresponding to the tuple
\[
((2 3), (2 3), (4 5), (4 5), (1 2), (1 2), \ldots, (1 2), (1 2)(3 4)(5
6), (1 2)(3 4)(5 6)) \, .
\]

We  let $Y_{g,n}$ be the closure in $X_{g,n}$ of points of the
form $C$, $\{a_1,a_2, \ldots, a_n\}$,\\
$\{b_1,b_2, \ldots, b_n\} $, where
$C$ is a cover of $\mathbb{P}^1$ corresponding to a point on
the chosen component of the Hurwitz scheme, and $\{a_1,a_2, \ldots, a_n\}$,
$\{b_1,b_2, \ldots, b_n\} $ are the inverse images of the 
two points of $\mathbb{P}^1$ over which the cover is not simply
ramified. Since the dimension of each component of the Hurwitz
scheme is $2g + 2n$ and $dim(Aut(\mathcal{P}^1)) = 3$, it follows that
$dim(Y_{g,n}) = 2g + 2n -3$.

To prove the statement about the  $\tau(Y_{g,n})$,
  we consider a degeneration of
the cover such that three of the simply ramified points come together
(generically) at a point. For instance, in the above example we would
let $p_5$, $p_6$ and $p_7$ come together.  If $g \geq 2$, the special
fibre of the stable model of the degenerating family of genus $g$
curves then consists of two smooth components, one of them a curve of
genus $g-1$ which is a cover of $\mathbb{P}^1$ of the same type, and
the other a generic elliptic curve. If $g=1$, we obtain a nonconstant
family of elliptic curves and so the statement is true in this case.
The statement for $g=2$ follows, since we may then assume that both
components are generic elliptic curves.

If $2 < g \leq 2n$, we may assume using induction that the point of
intersection of the two components is a generic point on the curve of
genus $g-1$ (also generic).
 Thus $dim(\tau(Y_{g,n})) \geq dim(\tau(Y_{g-1,n})) + 3$,
and so $\tau|_{Y_{g,n}}$ must be dominant. Finally, if $g > 2n$ we
see that $dim(\tau(Y_{g,n})) \geq dim(\tau(Y_{g-1,n})) + 2$.
Since  $dim(Y_{g,n}) = 2g + 2n -3$ and $dim(\tau(Y_{2n,n}))= 4n -3$,
it follows that we must
have equality for all $g \geq 2n$.

\end{proof}
%\subsection{}
The basic idea for the proof of indecomposability is now the same as
that used for the 4-configuration i.e. we specialize to a
hyperelliptic curve and then use the results of sections \ref{sec:hyp}
and \ref{sec:decomp}. In somewhat more detail, the argument is as
follows.  By Lemma \ref{lemma:spec} we may construct a family of
smooth curves $\mathcal{C} \to S$ of genus $g$ and a divisor
$\mathcal{D}= \sum_n A_i - \sum_n B_i$ on $\mathcal{C}$, with $S$ a
smooth curve with a distinguished point $s$, such that the fibre over
the generic point, $(C, D)$ corresponds to a generic point of $Y_{g,n}$
and the fibre over $s$, $(C',D')$, corresponds to the special points
of $Y_{g,n}$ in Lemma \ref{lemma:spec} (iii).

Let $W$ be the subspace of $CH^g(J(C'),1)$ obtained by mapping
$V(C,D)$ to $CH^g(Pic^{n+1}(C),1)_{\mathbb{Q}}$, specialization to
$CH^g(Pic^{n+1}(C'),1)_{\mathbb{Q}}$, and then translation by
$-(n+1)[w_1]$.  Note that specialization is always defined here, since
we are free to modify $f$ by a nonzero constant.

By the condition on $D'$ it follows that $f'$, the specialization of
$f$ must be a Weierstrass function on $C'$ corresponding to the
divisor $2w_1-2w_2$. The components of the support of elements of $W$
are all images of the various diagonals in ${C'}^{n+1}$, and the
function on each component is $f'$.  Lemma \ref{lemma:hyp} then allows
us to assume that modulo decomposable elements each element of $W$ is
a sum of translates of the basic hyperelliptic cycle $K$.  It is then
clear that the results of sections \ref{sec:hyp} and \ref{sec:decomp}
imply that if $g \geq n+ 1$, the general element of $V(C,D)$ (with
$(C,D)$ corresponding to a generic point of $Y_{g,n}$) is
indecomposable, provided that the following hypothesis is satisfied.
\begin{hyp}
  The image of $W$ in $CH^g_{\mathrm{ind}}(J(C'),1)_{\mathbb{Q}}$ is
  the vector space spanned by the cycle $K*([t_1 - w_1] - e)*([t_2 -
  w_1] - e)* \ldots *([t_{n-1} -w_1] - e)$.
\end{hyp}
\begin{rem}
The difficulty in verifying the hypothesis is purely combinatorial; everything
can be computed explicitly for any given $n$. We have verified it using a
computer for $2 \leq n \leq 6$. (In all these cases, the kernel
of the map from $W$ to $CH^g_{\mathrm{ind}}(J(C'),1)_{\mathbb{Q}}$
consists of those cycles which do not contain the (small) diagonal of
$C^{n+1}$ in their support.) Using the dimension formula from Lemma 
\ref{lemma:spec}, we see that $CH^g_{\mathrm{ind}}(J(C),1)_{\mathbb{Q}}$ 
is nonzero for
a generic curve with $3 \leq genus(C) \leq 12$. As in the case
of the $4$-configuration, Pontryagin product with zero cycles
can be used to prove that for these cases $CH^g_{\mathrm{ind}}(J(C),1)_{\mathbb{Q}}$
is in fact uncountable.
\end{rem}

\section{Higher analogues  in lower genus} \label{sec:higher}

In the first part of the section we construct elements $F \in
CH^3(J(C),4-g)$, where the genus $g \leq 2$. The cycles $F$ are
natural generalizations of the 4-configuration and we expect that they
should be strongly indecomposable in the sense that we explain below.
In the second part we study  a certain cycle $B \in CH^3(J(C),2)$ where $C$
is a bielliptic curve of genus $2$.  We show that $B$ is (weakly)
indecomposable using Lewis' criterion, which we show to hold by means 
of the same kind of proof as the one which was given in (1.3).

\subsection{Higher analogues of the 4-configuration in lower genus}

Bloch's groups $CH^p(X, n)$ can be described by means of chains built
from those integral subvarieties of codimension $p$ in $ X \times
(\mathbb P^1 - \{1\}) ^{n} $ which meet all the cubical faces over $0$
or over $\infty$ properly.  Consider a semistable degeneration $B'$ of
an abelian variety $B$ with trivial extension class, then $B' = A
\times \mathbb C^{*} $ with compactification $\bar B' = A \times \mathbb
P^1 $ and thus one may expect that $CH^{m}(A,n)$ should retain memory
of the properties of $CH^{m}(B, n-1)$.  This was the guess that
prompted our construction of the 4-configuration, which in a way can
be seen as being the memory of cycles studied in \cite{fakhruddin}.
In the same vein we describe now a unified procedure for building a
series of cycles $F$ in $CH^3(J(C),4-g)$, $g=genus(C)$ $ \leq 3$.  The
first is the 4-configuration and each element is the memory of the
preceding one.

Choose a class $ \epsilon $ of order $2$, consider a divisor $ D:=
(a_1 + a_2) - (b_1 + b_2 )$ on $C$ with $class (D) = \epsilon $, and
let $f$ be a rational function with $ div (f) = 2 D $. We embed $C \to
J(C)$ using the maps $\alpha _1 (x) = [x-a_1]$, $\alpha _2 (x) = [
-x+a_2]$, $\alpha _3 (x) = [x-a_1]+ \epsilon$, $\alpha _4 (x) =
[-x+a_2]+ \epsilon$, and define $C_i := \alpha _i (C)$. The useful
property here is the fact that on $C_1$ the point which comes from
$a_i$ coincides with the point on $C_2$ from $a_{ i\pm 1} $. The same
happens for $b_i$ with respect to the curves $C_1$ and $C_4$, and
things stay the same for every other curve instead of $C_1$.
        
We start from genus $3$, here we take $\beta _i = C \to J(C) \times
\mathbb P^1 $ to be the map $ x \to (\alpha _i ( x), f(x)^{s(i)}) $,
where $s(i)$ means $\mp 1$, according to the parity of i.

The curves $K_i := \beta _i(C) $ in $J(C) \times \mathbb P^1 $ meet
properly the cubical faces, $\sum_{ i = 1} ^{4} K_i$ is indeed a cycle
and this is our element
\[
F(f): =\sum_{ i =1} ^{4} K_i \in CH^3(J(C),1) \, .
\]

\noindent Clearly $ F(f)$ is equivalent to the  the 4-configuration
as it was previously defined.

To go higher in Chow groups we need at each step a new and convenient
rational function. It can be found by imposing restrictions on $D$ by
means of conditions on $f$ and this is the reason why we need to go
down in genus.

In genus $2$ we require $ \epsilon = [ w_1 - w_2] $, the distinguished
ramification points of the Weierstrass double cover $ h : C \to
\mathbb P^1$.  The condition is
\begin{equation}
 f(w_1) = f(w_2)  \tag{+} 
\end{equation}
which is satisfied by a 1- dimensional family of $a$'s and $b$'s as
above. Weil reciprocity and (+) yield
\begin{equation}
 (h(a_1) h(a_2))^2 = (h(b_1) h(b_2))^ 2 \, .  \tag{*}  
\end{equation}
  
We may change $ h $ by a multiplicative constant so as to have
\begin{equation}
1= \ (h(a_1) h(a_2))^2 = (h(b_1) h(b_2))^ 2  \, . \tag{**}
\end{equation}
Writing now $ \beta _i ( x) := (\alpha _i ( x), f(x)^{s(i)},
h(x)^{2s(i)}) $ then $K_i := \beta _i (C)$ is a curve in $J(C)\times
\mathbb P^1 \times \mathbb P^1 $ which meets properly the cubical
faces.  It is easy to check that
\[
F(f,h) := \sum_{ i = 1} ^{4} K_i \in CH^3(J(C),2) \, ,
\]
for instance one can see that the boundary is trivial over the point
of $C_1$ which comes from $w_1$ by realizing that it coincides with
the point on $C_3$ from $w_2$, and by using then the condition $f(w_1)
= f(w_2)$.

\medskip

In genus $1$ having fixed the divisor $D$ of class $\epsilon$ as
before we choose further two rational functions $ h_1 $ and $ h_2 $ of
degree $2$ on $E$ both ramified over $0$ and $\infty$.  Let $ div
(h_i) = 2(q'_{i} -q''_{i})$.  We require our choice to satisfy:
\begin{align}
  & [q'_{1}  -q''_{1}]  = [q'_{2}  -q''_{2}] =  \epsilon \tag{1} \\
  & h_{1} (q'_{ 2}) = h_{1} (q''_{ 2}) \quad and \quad h_{2} (q'_{ 1})
  = h_{2} (q''_{ 1})  \, .\tag{2}
\end{align}

This can be done, see \cite{collino-jag}. Here the conditions on $f$
read
\begin{equation}
f(q'_i) = f(q''_i)  \quad   i =  1 \,  , \, 2 \, . \tag{+} 
\end{equation}
As it was before this implies
\begin{equation}
 (( h_i (a_1)   h_i (a_2) )^2 = ( h_i (b_1)   h_i(b_2) )^2  \, .  \tag{*}
\end{equation}
We normalize the choice of the rational functions $h_i$ by the
request:
\begin{equation}
 1= (( h_i (a_1)   h_i (a_2) )^2 = ( h_i (b_1)   h_i(b_2) )^2 \, . \tag{**}
\end{equation}
\noindent Using  
$ \beta _i ( x) := ( \alpha _i ( x), f(x)^{s(i)}, h_1 (x)^{2s(i)}, h_2
(x)^{2s(i)} ) $, we obtain
\[
F(f,h_1,h_2) := \sum_{ i = 1} ^{4} K_i \in CH^3(E, 3) \, .
\]
  
The 4-configuration $ F(x,f,h_1,h_2) \in CH^3(\mathbb C, 4)$ is
constructed in the same way as $ F(f,h_1,h_2)$ was.  One thinks of
$\mathbb P^1$ as having the structure of a degenerate Jacobian, given
by the choice of a standard parameter $x$. The opposite map from
$J(C)$ to $J(C)$ becomes $ x \to x ^{-1} $ and the sum of points
corresponds to product of the coordinates. Translation by $ \epsilon $
is here multiplication by $-1$.  The condition $class (D) = \epsilon $
reads $a_1 a_2 = - b_1 b_2 $ and $f := (x- a_1)^2 (x - a_2)^2 (x-
b_1)^{-2} (x - b_2)^{-2}$.  We may take explicitly: $h_1 = c_{1} (x-
1)^2 (x +1 )^{-2}$ , $h_2 = c_{2} (x- i)^2 (x +i )^{-2}$. The
requirements on $f$ are

\begin {equation}
  f(1) = f(-1) \, , \quad f(i) = f(-i) \tag{+} \,
\end{equation}
Our choice yields also $: \, f(0) = f(\infty) $.  One has
\begin{equation}
( h_i (a_1) h_i(a_2))^2 = (h_i (b_1) h_i (b_2))^2 \quad  , i = 1, 2. \tag{*}
\end{equation}

Choose and fix the constants $c_{1}$ and $c_{2}$ so that it is
\begin{equation}
1= ( h_i (a_1) h_i(a_2))^2 = (h_i (b_1) h_i (b_2))^2 \, , i = 1, 2. \tag{**}
\end{equation}

With this dictionary the maps $\alpha_i: \mathbb P^1 \to \mathbb P^1$
are defined as before (for instance $\alpha_1 (t) = ta_1^{-1}$ ), here
the range of $\alpha_i$ should be understood as the replacement of
$J(C)$. In this way the maps $ \beta : \mathbb P^1 \to (\mathbb P^ 1)
^{4} $ are here $ \beta _i ( x)= ( \alpha _i ( x), f(x)^{s(i)}, h_1
(x)^{2s(i)}, h_2 (x)^{2s(i)} ) $. Setting again $K_i := \beta_i
(\mathbb P^1)$ our cycle is then 
\[
F(x,f,h_1,h_2) := \sum_{ i = 1}^{4} K_i \in CH^3(\mathbb C^{*}, 4) \, .
\]

\begin{rem} One may define an element in $CH^a(X,b)$
  to be strongly indecomposable if it is not in the image of $CH^{a-1}
  (X,b-1) \otimes \mathbb C^ {\times}$.  We think that the higher
  4-configurations are good candidates to strong indecomposability.
\end{rem}

\subsection{The B-configuration}

Following the terminology of \cite{lewis-crelle} we define the group
of (weakly) decomposable cycles in $ CH^{p}(X,2) \simeq
H^{p-2}(X,\mathcal K_{p})$ to be :
\[
CH_{\mathrm{dec}}^p(X,2) := Im \big\{ K_2(\mathbb C) {\otimes}
CH^{p-2}(X) \longrightarrow CH^p(X,2)\big\} \, ,
\]
and thus the{{\it} indecomposable\/} group is
\[
CH_{\mathrm{ind}}^k(X,2) :=  CH^{k}(X,2)\big/
CH_{\mathrm{dec}}^k(X,2) \, .
\] 

It is known that translations on an elliptic curve $E$ act trivially
on $CH_{\mathrm{ind}}^2(E,2)$, see \cite[3.10]{goncharov-levin}.  We
show that on the contrary translations on a genus $2$ bielliptic
Jacobian $J(G)$ operate non trivially on $CH_{\mathrm{ind}}^3(J(G),2)$.
Our procedure is similar to the one which we have applied above  in (1.3).
We deal here with the { \it
  B-configuration } which is shown to be indecomposable by checking
Lewis' condition on a cycle ${\mathfrak B }$ of $CH^{3}(J(G)\times
G,2)$.

S. Bloch in his seminal memoir \cite{bloch-lac} constructed certain
elements $S_b \in \Gamma(E,\mathcal K_{2}))$ associated with a point
$b$ of finite order on an elliptic curve $E$.  He proved that the real
regulator image of $S_b$ is not trivial for some curves with complex
multiplication, and thus it is not trivial in general.

Consider a bielliptic curve $G$ of genus $2$ with associated map
$\delta_G :G \to E_1 $, and let $a: G \to J(G)$ be the Abel Jacobi
map.  In this way $Z(b) := a_{\ast} \delta_G^{\ast} (S_b)$ is a cycle
in $CH^{3}(J(G),2)$.  Translation by an element $t \in Pic^0 (G)$ maps
$Z(b)$ to the cycle $Z_t(b)$, our aim is to prove
\begin{thm}
  The B-configuration $B(t) := Z_t(b) - Z(b)$ is indecomposable for
  generic $t$.
\end{thm}
\noindent Note that $B(t)$ has trivial regulator image.
\begin{proof}
  Consider the cycle $G \times Z(b) \in CH^{2}(G \times G,2)$.  The
  straight embedding $\sigma := id \times a: G \times G \to G \times
  J(G)$ maps it to ${\mathfrak S } := \sigma_{\ast} (G\times Z(b))$ in
  $CH^{3}(G \times J(G),2)$.  The twisted embedding $\tau (t,x) :=
  (t,a(x)+(t-w_1)))$ gives instead ${\mathfrak T }:= \tau_{\ast}
  (G\times Z(b)) $, with section ${\mathfrak T }_{\ast}(t) = Z_t(b)$,
  and therefore $B(t)$ is the section at $t \in G$ of ${\mathfrak B
    }:= {\mathfrak T } - {\mathfrak S }$.
  
  We use the same type of notations as we did in part \ref{subsec:zm}, in
  particular the holomorphic form ${\omega}_i^{J}$ comes from $E_i$.
  We need to consider also the forms $\nu
  :=\bar{\omega}_1^{J}\wedge{\omega}_2^{J} $ on $J(G)$ and
  $\bar{\zeta_2}$ on $G$.  The procedure of \ref{subsec:zm} gives here again:
  (i) $<R({\mathfrak B }),\bar{\zeta_2} \wedge \nu > \neq 0$.  The
  Neron-Severi space of divisors with rational coefficients on $J(G)$
  is isomorphic to the same space on the product of the two associated
  elliptic curves.  On the gneric bielliptic Jacobian $\nu$ is
  orthogonal to the Neron-Severi group, because it is orthogonal to
  the elliptic curves.  A look at the proof of the main theorem of
  \cite{lewis-crelle} shows that this property of $\nu$ and (i) imply that the
  generic section ${\mathfrak B }_{\ast}(t)$ is indeed weakly
  indecomposable.
\end{proof}

\smallskip
\begin{ack}
The second author would like to thank R.~Sreekantan for some
useful correspondence.
\end{ack}

\end{document}